\documentclass[a4paper, 11pt]{article}

\usepackage[margin=1.2in]{geometry}
\usepackage{parskip}
\usepackage[titletoc,title]{appendix}

\usepackage{array}
\usepackage{caption}
\usepackage{rotating}
\usepackage[breaklinks, hidelinks, colorlinks = true, linkcolor = blue, urlcolor = blue, anchorcolor = blue, citecolor = blue]{hyperref}
\usepackage{blkarray, bigstrut}
\usepackage{multirow}
\usepackage{booktabs}
\usepackage{siunitx}
\usepackage{mathtools}
\usepackage{amsmath}
\usepackage{amssymb}
\usepackage{amsthm}
\usepackage{complexity}
\usepackage{natbib}
\usepackage{graphicx}
\usepackage{pgf}
\usepackage{pgfplots}
\usepgfplotslibrary{groupplots}
\usepackage{tikz}
\usetikzlibrary{matrix,patterns,positioning,decorations,arrows,shapes,decorations.markings,calc}
\usepackage[latin1]{inputenc}
\usepackage{verbatim}
\usepackage{subfiles}
\usepackage{placeins}
\usepackage{pgfplots}
\usepackage{color}
\usepackage{float}
\usepackage{enumerate}
\usepackage{emptypage}
\usepackage{subcaption}
\usepackage{setspace}
\usepackage[export]{adjustbox}
\usepackage{subfiles}
\usepackage{algorithm}
\usepackage{algpseudocode}
\usepackage{todonotes}

\newcolumntype{H}{>{\setbox0=\hbox\bgroup}c<{\egroup}@{}}

\makeatletter
\pgfplotsset{
    groupplot xlabel/.initial={},
    every groupplot x label/.style={
        at={($({\pgfplots@group@name\space c1r\pgfplots@group@rows.west}|-{\pgfplots@group@name\space c1r\pgfplots@group@rows.outer south})!0.5!({\pgfplots@group@name\space c\pgfplots@group@columns r\pgfplots@group@rows.east}|-{\pgfplots@group@name\space c\pgfplots@group@columns r\pgfplots@group@rows.outer south})$)},
        anchor=north,
    },
    groupplot ylabel/.initial={},
    every groupplot y label/.style={
            rotate=90,
        at={($({\pgfplots@group@name\space c1r1.north}-|{\pgfplots@group@name\space c1r1.outer
west})!0.5!({\pgfplots@group@name\space c1r\pgfplots@group@rows.south}-|{\pgfplots@group@name\space c1r\pgfplots@group@rows.outer west})$)},
        anchor=south
    },
    execute at end groupplot/.code={%
      \node [/pgfplots/every groupplot x label]
{\pgfkeysvalueof{/pgfplots/groupplot xlabel}};  
      \node [/pgfplots/every groupplot y label] 
{\pgfkeysvalueof{/pgfplots/groupplot ylabel}};  
    }
}

\def\endpgfplots@environment@groupplot{%
    \endpgfplots@environment@opt%
    \pgfkeys{/pgfplots/execute at end groupplot}%
    \endgroup%
}
\makeatother

\usepackage{pdflscape}

\usepackage[normalem]{ulem}

\title{Benders Decomposition for Robust Tactical Railway Crew Scheduling}

\author{ B.T.C. van Rossum$^{1, \star}$, T. Dollevoet$^1$, D. Huisman$^{1,2}$
\vspace{0.1cm}\\
\small{$^1$Econometric Institute and Erasmus Center for Optimization in Public Transport}\\ 
\small{Erasmus University Rotterdam, The Netherlands}
\vspace{0.1cm}\\
\small{$^2$Digitalization of Operations, Netherlands Railways} \\
\small{Utrecht, The Netherlands}
\vspace{0.1cm}\\
\small{$^\star$Corresponding author} \\
\vspace{0.1cm} \\
\small{vanrossum@ese.eur.nl, dollevoet@ese.eur.nl, huisman@ese.eur.nl}
}

\date{October 8, 2024}

\begin{document}

\tikzset{->-/.style={decoration={
  markings,
  mark=at position #1 with {\arrow{>}}},postaction={decorate}}}

\tikzstyle{block} = [rectangle, draw, 
 text centered, rounded corners, inner sep = 0.2cm]

\tikzset{->-/.style={decoration={
  markings,
  mark=at position #1 with {\arrow{>}}},postaction={decorate}}}
  
\maketitle

\begin{abstract}
\noindent 
We consider robust tactical crew scheduling for a large passenger railway operator, who aims to inform crew early on about their work schedules while also maintaining the ability to respond to changes in the daily timetables. To resolve this conflict, the operator considers a template-based planning process, templates being time windows during which duties can later be scheduled. The goal is to select a cost-efficient set of templates that is robust with respect to uncertainty in the work to be performed in the operational phase. A set of templates is deemed robust when few excess duties are required to cover all work in the operational planning phase. To enable the construction of efficient template-based rosters, we impose several template rostering constraints that proxy the actual rostering rules of later planning steps. We propose a two-phase accelerated Benders decomposition algorithm that can incorporate these restrictions. Computational experiments on real-life instances from Netherlands Railways, featuring up to 948 tasks per day, show that historical planning information can be used to obtain robust templates and that parsimonious solutions can be obtained at negligible extra costs. Compared to a literature benchmark, our Benders decomposition method solves three times as many instances without rostering constraints to optimality. 
\vspace{5mm}
\end{abstract}

\section{Introduction}
\label{sec:intro}

Passenger railway operators engaging in crew planning must balance two conflicting goals. On the one hand, they wish to provide certainty to crew members about their working schedules early on. On the other hand, they like to maintain the flexibility to modify their planning in response to unforeseen changes on the other hand. In particular, crew members wish to be informed about their rosters already in the tactical planning phase, months before the day of operation, in order to schedule their personal lives around their work. At the same time, the exact tasks, i..e, scheduled train trips, to be performed are not known until the operational phase, at most several weeks before the day of operation. While a yearly timetable and rolling stock schedule are available, they remain subject to changes resulting from, for example, scheduled maintenance works and events. 


To handle these conflicting objectives, Netherlands Railways (NS), the largest passenger railway operator of the Netherlands, is considering a template-based planning process. Here, a template is a time window specifying the hours during which a crew member can be called upon to perform a duty, i.e., day of work. The proposed process starts with the construction of a template-based roster in the tactical planning phase. The template-based roster is communicated to crew several months before the day of operation and is considered fixed from that point onwards. The roster is detailed to a duty-based roster in the operational planning phase, once the timetable and rolling stock schedule for specific days are finalised. When the work to be performed is known, duties are generated that cover this work and can be performed by crew members within the time windows specified in their roster. Costly excess duties are available in case the template-based rosters fail to provide sufficient capacity.

The template-based planning process has several advantages. First, templates provide the operator a way of accounting for the uncertainty regarding the work to be performed. Since the length of template time windows typically exceeds the length of an average duty, a single template can accommodate multiple duties in different realisations of the daily timetable. Second, template-based rosters offer strong guarantees to crew members. They are able to schedule their personal lives well in advance, knowing that they will not be called upon to work outside the hours specified in the templates. Finally, there is no need for many costly (re-)planning steps, since duties are only generated once in the operational planning phase.

Motivated by NS, we study the viability of such a template-based process and consider the tactical problem of choosing templates that are both cost-efficient and robust with regards to uncertainty in the work to be performed in the operational phase. To this end, we use the concept of recoverable robustness, and say that a set of templates is robust when few or no excess duties are required to cover all work in the operational planning phase. We assume that some historical planning information is available to inform the selection of robust solutions. To ensure that efficient template-based rosters can be constructed in subsequent planning steps, we also enforce several template rostering constraints that proxy the actual rostering rules. For example, we limit the number of unique template types, reserve templates, and early or late templates. While we focus on a large passenger railway operator, the setting where preliminary capacity plannings are constructed before detailed work assignments is encountered in many (crew) scheduling applications.

While there is a growing interest in robust railway planning, few authors study robust railway crew planning \citep{lusby2018survey}. In addition, most authors focus on generating duties that are robust to daily disturbances, whereas our goal is to construct a capacity plan that is robust with respect to minor changes in the timetable. Our work is most closely related to that of \cite{rahlmann2021robust}, who study robust tactical crew planning for a rail freight operator. Using a column generation heuristic, they choose robust templates based on a finite set of historical planning scenarios. We build on their work by imposing a wide range of template rostering constraints, necessitating a novel solution approach that is able to incorporate these restrictions. Moreover, our passenger railway application features a higher number of tasks per day and more tasks in a typical duty.

This paper has four main contributions. First, we extend the problem setting of 
\cite{rahlmann2021robust} by imposing template rostering constraints. Second, to handle these restrictions, we propose an accelerated two-phase Benders decomposition algorithm. Third, we apply our method to real-life instances from NS featuring up to 948 tasks per day. We find that historical planning information can be used to obtain robust crew schedules and that parsimonious crew plans with a limited number of unique template types can be obtained at negligible extra costs. Fourth, we benchmark our approach against the column generation heuristic of \cite{rahlmann2021robust}. The Benders decomposition algorithm solves three times as many large instances without rostering constraints to optimality, at the price of a two-fold increase in computing time. 

The remainder of this paper is structured as follows. We introduce the robust tactical crew scheduling problem in Section~\ref{sec:problem} and provide an overview of related literature in Section~\ref{sec:lit}. We motivate and present a mathematical formulations in Section~\ref{sec:recoverable}, and propose an accelerated two-phase Benders decomposition algorithm in Section~\ref{sec:benders}. In Section~\ref{sec:benchmark}, we briefly present the benchmark algorithm of \cite{rahlmann2021robust}. We conduct extensive computational experiments on instances from NS, discussing computational performance and practical insights in Sections~\ref{sec:computational} and \ref{sec:practical}, respectively. We conclude in Section~\ref{sec:conclusion}.

\section{Problem Description}
\label{sec:problem}

The goal of crew planning is to assign all available work to crew members in an efficient manner while respecting all relevant labour rules. The work to be performed consists of tasks, i.e., scheduled train trips. In crew scheduling, tasks are aggregated to duties, i.e., feasible sequences of tasks constituting days of work. In crew rostering, duties are assigned to crew members in the form of rosters, detailing which duties are to be performed by which employee at each day of the planning period. Typically, crew scheduling and crew rostering are performed sequentially and in the tactical planning phase. Later, in the operational planning phase, duties and rosters are modified to account for changes in the timetable. The current crew planning process at NS follows this structure. The detailed planning steps in the tactical phase limit the flexibility of the operator in the operational phase, while the many roster updates provide little certainty to crew members. We refer to \cite{abbink2018crew} and \Citet{rossum2024railway} for more information on the planning process at NS.

\begin{figure}[H]
\centering
\begin{adjustbox}{max width= \textwidth}
\begin{tikzpicture}
	\draw (2.5, 4.5) node {\textbf{Tactical}};
	\draw (7, 4.5) node {\textbf{Operational}};

	\draw[thick, dotted] (5, 0.5) -- (5, 4);

	\tikzset{every node/.style={draw, rectangle, node distance=1cm, minimum height=2em, inner sep=4pt, anchor=south west}}



	\draw (3, 1) node (6) {Rosters};
	\draw node[left = of 6] (5) {Templates};
	\draw node[above = of 5, fill=black!15] (0) {Historical data};
	\draw node[right = of 6] (7) {Duties};
	\draw node[above = of 7, fill=black!15] (8) {Tasks};

	\path[->,>=stealth'] (0) edge (5); 
	\path[->,>=stealth'] (5) edge (6);      
	\path[->,>=stealth'] (6) edge (7);
	\path[->, >=stealth'] (8) edge (7);

\end{tikzpicture}
\end{adjustbox}
\caption{Proposed crew planning process. Crew planning decisions are indicated by white rectangles, inputs by shaded rectangles.}
\label{fig:process}
\end{figure}
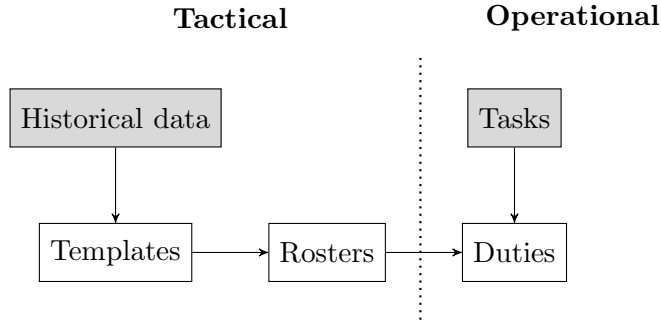

The proposed planning process, illustrated in Figure~\ref{fig:process}, effectively reverses the order of the two planning steps. A template-based crew rostering step takes place in the tactical planning phase, and all crew scheduling is delegated to the operational planning phase. In particular, in the tactical planning phase we do not schedule duties but build a roster based on templates only. Here, a template is a combination of crew base and time window, representing a placeholder for a duty that is to be performed by a crew member of this base during the specified time window. In other words, a template specifies an interval during which a crew member will be called upon to work. The template-based roster is communicated to the crew before the start of the operational planning phase. The operator may not deviate from the working hours communicated here throughout the remainder of the planning period, enabling crew members to better schedule their personal lives. During the course of the operational phase, the daily timetables are finalised and the tasks to be performed on each day of the planning period are revealed. Duties must be generated that cover these tasks and can be performed with the capacity provided by the templates in the roster. If a duty cannot be performed by a crew member during his scheduled shift, as specified by the template, the duty becomes an excess duty that will be executed by an extra crew member at a high penalty cost. 

We focus on the problem encountered in the tactical planning phase, in which several objectives must be balanced. First, the railway operator desires an efficient crew plan with low labour costs, and wishes to minimise the number of selected templates. Second, it must be possible to construct balanced rosters with the chosen templates, and hence some rostering compatibility constraints must be taken into account. Third, the templates must be chosen such that the number of excess duties required in the operational phase is minimal. In other words, the templates must be robust against many possible realisations of the tasks to be performed in the operational phase. We assume that historical planning information is available in the tactical planning phase, which can be leveraged to select robust templates. In short, the goal of the robust tactical crew scheduling problem is to select templates that satisfy roster compatibility constraints and minimise a weighted sum of template and excess duty costs. 

We discuss the restrictions encountered in the tactical and operational planning phases in more detail in Sections~\ref{subsec:problem_tactical} and \ref{subsec:problem_operational}, respectively. 

\subsection{Tactical Phase}
\label{subsec:problem_tactical}

In the tactical planning phase, we are provided with a pre-specified list of template types and must determine how often to use each template type. Each template acts as a placeholder for duties, and is characterised by a crew base, day of the week, earliest start time, and latest end time. The list contains regular templates of a fixed length, as well as a reserve template spanning the entire day. Based on the chosen templates, each crew base will construct a preliminary roster that is communicated to crew members before the start of the operational phase. We focus on the problem of choosing templates, and consider the actual crew rostering step to be outside the scope of this paper. 

To ensure that efficient rosters can be constructed, the chosen templates at each crew base must satisfy several constraints. First, we impose some structure within the roster, increasing the opportunities for exchanging shifts, by limiting the number of unique template types that can be selected. Second, we limit the percentage of templates starting early (before 6:00h) or ending late (after 24:00h). The maximum percentage of early and late templates depends on the day of the week, as the percentage of night duties required is generally higher in the week-end. Third, the percentage of reserve templates may not exceed a certain upper bound. Finally, one could impose a maximum capacity constraint if the number of available crew members is known to be limited. 

\subsection{Operational Phase}
\label{subsec:problem_operational}

Throughout the operational phase, the timetable and rolling stock schedule for each day of the planning period undergo changes to accommodate for, among other things, maintenance works and events. Once the schedules for a particular day are finalised, the set of tasks to be performed on that day is known. All tasks must be covered by duties, i.e., days of work. Each duty can either be assigned to a template or turned into an excess duty. When a duty is assigned to a template, it is performed by a crew member scheduled to work in a time window spanning the duty, and no additional cost is incurred. A duty that cannot be assigned to a template becomes an excess duty, to be performed by an excess crew member at a high penalty cost. The goal of the operational planning phase is to construct, for each day of the planning period, a set of duties that cover all tasks and minimise the sum of penalty costs and a variable cost depending on the total workload. 

Each duty is a sequence of tasks that can be performed consecutively and must satisfy several duty rules. Two tasks can be performed consecutively when they end and start at the same station, respectively, and allow for a sufficient transition time in case they take place on different rolling stock units. The duty rules are as follows. First, each duty starts and ends at the same crew base. Second, the total length of a duty, defined as its ending minus starting time, may not exceed a maximum length that depends on the starting time of the duty. Third, each duty must include a meal break of at least 30 minutes at a station with a dedicated canteen. To position the meal break approximately halfway the duty, no duty may contain a period longer than 5.5 hours without a meal break.

\section{Literature Review}
\label{sec:lit}

This work belongs to a large body of literature on crew planning in public transport. Crew planning is generally decomposed into crew scheduling (or, crew pairing in airline applications) and crew rostering. We refer the interested reader to \cite{heil2020railway} for a recent review on railway crew scheduling, and to \cite{caprara1998modeling} and \cite{kohl2004airline} for introductions to railway and airline crew rostering, respectively. \cite{abbink2018crew} provide an overview of crew planning at NS.

A considerable stream of research focuses on template-based scheduling approaches in tactical railway crew planning. \cite{borndorfer2017integration} integrate duty scheduling and rostering, coupling the two problems through duty templates. Similar to our work, \cite{rahlmann2021robust} consider the goal of selecting robust templates in tactical crew planning for a freight railway operator. \cite{breugem2022column} study crew replanning in case of scheduled disruptions. The use of templates ensures that replanned rosters do not deviate much from the originally scheduled working hours. Finally, \Citet{rossum2024railway} consider a planning process that is highly similar to ours. They focus on the operational planning phase, where individual template-based rosters are fixed, and aim to attain temporal individual fairness of the scheduled duties. 

The large impact of disruptions on operational performance of railway systems has led to a large body of research on robustness in railway planning. We refer to \cite{lusby2018survey} and \cite{bevsinovic2020resilience} for recent reviews on robust planning. Within crew planning, robustness is less well-studied, with most authors focusing on constructing duties that are robust to daily disruptions. \cite{veelenturf2016quasi} consider real-time crew rescheduling in case of disruptions. They introduce the notion of `quasi-robust' rescheduled duties: duties that are feasible in the most optimistic duration of the disruption, but can be converted to feasibe duties in case the disruption lasts longer. \cite{fuentes2019hybrid} argue that robustness can be incorporated by increasing the minimum transfer time between tasks, or by penalising deadheading movements and short transfers. Since we focus on the tactical planning phase where no duties are scheduled, robustness with regards to daily disruption is outside the scope of our work. 

A larger stream of research considers robust airline crew pairing. \cite{antunes2019robust} study airline crew pairing with the goal of generating pairings that are less likely to be disrupted by small delays. \cite{muter2013solving} consider a setting where extra flights may be scheduled during operation, and aim to generate pairings that can accommodate extra flights with minimal disruptions. They propose a row-and-column-generation algorithm to solve the problem. \cite{weide2010iterative} provide an iterative solution approach for integrated airline crew pairing and aircraft routing, aiming to generate solutions that are robust against operational delays. \cite{dunbar2014integrated} further integrate re-timing decisions and aim to minimise the effect of delay propagations. Finally, \cite{shebalov2006robust} build robust crew pairings by maximising the number of so-called `move-up crew', crew members whose tasks can be freely swapped in operation. They solve this problem using a combination of Lagrangian relaxation and column generation. Since we consider robustness against timetable modifications, instead of robustness against operational delays, our work is most similar in spirit to that of \cite{muter2013solving}.

Our work is most closely related to that of \cite{rahlmann2021robust}, who consider robust tactical crew planning for a railway freight operator. Similar to us, they aim to choose a set of templates that are able to fit efficient duties in multiple demand scenarios. They propose a column generation approach with a two-layered pricing structure to solve the problem. Whereas we assume a pre-specified list of template types to choose from, they allow templates to start at arbitrary times of the day. In addition to their work, we impose template rostering constraints to incorporate some aspects of crew rostering, necessitating a novel solution approach.

\section{Recoverable Robust Formulation}
\label{sec:recoverable}

Two modelling choices underlie our mathematical formulation. First, we use the concept of recoverable robustness to model the interaction between the tactical phase and operational phase. Recoverable robustness, first introduced by \cite{liebchen2009recoverable}, pertains to two-stage robust optimisation problems, where first-stage decisions might be rendered infeasible by a realisation of some uncertain data, but feasibility can be recovered at minor cost through second-stage decisions. This perfectly captures the relation between the tactical and operational phase: a set of templates is deemed recoverable robust when few excess duties are required to perform all tasks in the operational phase. The number of excess duties required in the second-stage can be determined by solving a (capacitated) crew scheduling problem. This approach suggests constructing a model that encapsulates both first-stage, tactical decisions, as well as second-stage, operational decisions.

Second, we capture the uncertainty regarding tasks to be performed in the operational phase using a finite set of scenarios. From one planning week to the next, tasks can be modified, added to, or removed from the schedule. Since detailed temporal and geographical information is needed to determine whether tasks can be contained in the same duty, such dynamics are hard to capture in a compact manner with traditional, polyhedral uncertainty sets. Instead, we use a scenario-based approach. The key assumption here is that a finite set of task scenarios, each corresponding to a possible realisation of tasks in the operational phase, is representative of actual operations. In other words, we assume that if a tactical solution performs well in all the considered scenarios, it is likely to perform well in future operations. This approach is viable in practice, since railway operators typically have access to a large set of historical planning information that can be used to construct a proper scenario set.

We finish with a brief discussion on how to construct scenario sets. While a large amount of historical planning information is typically available, attention must be paid to ensure that the chosen scenarios are both representative of future operations and contain sufficient variation to yield robust solutions. For example, historical data might be corrupted in case of replanning due to one-off events like scheduled maintenance works. In a scenario set, however, one rather aims to capture small-scale daily disturbances that are likely to be repeated, like additional shunting movements early in the morning or late in the evening. Finally, there is the question of how many distinct scenarios are required to obtain robust solutions. We will answer this question through extensive computational experiments in Section~\ref{sec:practical}.

\subsection{Mathematical Formulation}
\label{subsec:recoverable_formulation}

We now show how to formalise the modelling choices and provide a mathematical formulation of the robust tactical crew scheduling problem. We introduce the following notation. Let $P$ denote the set of templates and $c_p$ the cost of using template $p$. Denote the maximum number of chosen template types by $\Gamma$, and the set of template rostering constraints by $M$. All constraints we consider, such as the maximum percentage of early templates, can be modelled in a linear fashion. For each constraint $m \in M$, $E_m$ is a right-hand side constant, and coefficient $e^{m}_{p}$ denotes the left-hand side contribution of template $p$. Let $S$ denote the set of scenarios, and $K_s$ denote the set of tasks to be covered in scenario $s$. Let $\Delta^{s}_{p}$ denote the set of duties that can be assigned to template $p$ in scenario $s$, and let binary parameter $a_{pk}^{\delta s}$ indicate whether task $k \in K_s$ is covered by duty $\delta \in \Delta_{p}^{s}$. Finally, let $c_E$ denote the cost of an excess duty. 

Introduce binary and integer decision variables $z_p$ and $y_p$ to denote if and how often template $p$ is used, respectively. Moreover, let continuous variable $v^{s}_{p}$ denote the number of excess duties of template $p$ used in scenario $s$, and let $\eta$ denote the maximum recovery cost across all scenarios. Finally, let binary decision variable $x^{\delta s}_{p}$ denote if duty $\delta$ is assigned to template $p$ in scenario $s$. The recoverable-robust formulation of the tactical crew scheduling problem then reads:
\begin{subequations}
\begin{align}
\min \quad & \sum_{p \in P}  c_{p} y_{p} + \eta \label{eq:obj} \\
\text{s.t.} \quad & y_{p} \leq C z_{p} && \forall p \in P \label{eq:template_activation}  \\ 
& \sum_{p \in P} z_{p} \leq \Gamma \label{eq:template_types} \\
& \sum_{p \in P} e_{p}^{m} y_p \leq E_m && \forall m \in M \label{eq:template_mix}  \\ 
& \sum_{p \in P} c_E v^{s}_{p} \leq \eta && \forall s \in S \label{eq:recovery} \\ 
& \sum_{\delta \in \Delta_{p}} x_{p}^{\delta s} \leq y_{p} + v_{p}^{s} && \forall s \in S, p \in P \label{eq:excess} \\
& \sum_{p \in P} \sum_{\delta \in \Delta^{s}_{p}} a_{pk}^{\delta s} x_{p}^{\delta s} \geq 1 && \forall s \in S, k \in K_s \label{eq:cover} \\
& v_{p}^{s} \geq 0 && \forall s \in S, p \in P \label{eq:domain_excess} \\
& x_{p}^{\delta s} \in \mathbb{B} && \forall s \in S,  p \in P, \delta \in \Delta_{p} \label{eq:domain_duty} \\ 
& y_{p} \in \mathbb{N} && \forall p \in P \label{eq:domain_template_number} \\
& z_{p} \in \mathbb{B} && \forall p \in P. \label{eq:domain_template_selection}
\end{align}
\label{eq:model}%
\end{subequations}
The objective~(\ref{eq:obj}) is to minimise the sum of template costs and the maximum recovery cost. Constraints~(\ref{eq:template_activation})-(\ref{eq:template_types}), where parameter $C$ is an upper bound on the number of required templates, ensure that no more than $\Gamma$ different template types are used. Constraints~(\ref{eq:template_mix}) enforce the template rostering conditions, and constraints~(\ref{eq:recovery}) correctly set the maximum recovery cost. Constraints~(\ref{eq:excess}) ensure that duties of a template are covered by either available templates or excess duties, while constraints~(\ref{eq:cover}) ensure that each task is covered by at least one duty. Constraints~(\ref{eq:domain_excess})-(\ref{eq:domain_template_selection}) model the variable domains.

Model~(\ref{eq:model}) is challenging to solve for multiple reasons. First, it contains many integer and binary template variables, coupled to each other through weak big-$M$ type constraints~(\ref{eq:template_activation}). Second, the number of possible duty variables is already large for a single scenario, and this model contains multiple scenarios. Our attempts at column generation-based heuristics have proven unsuccessful, which is why we propose a Benders decomposition algorithm in Section~\ref{sec:benders} instead.

\subsection{Evaluation in Operational Phase}
\label{subsec:evaluation}

The solution quality of templates selected in the tactical planning phase is evaluated in the operational planning phase. On each day of the planning period, we aim to select duties that minimise the sum of excess duty penalties and variable workload costs. Let $K$ be the set of tasks to be covered on a particular day and $\bar{y}_p$ the number of templates of type $p$ selected in the tactical phase. The operational crew scheduling problem can then be formulated as follows:
\begin{subequations}
\begin{align}
\min \quad & \sum_{p \in P} \sum_{\delta \in \Delta_{p}} l_{p}^{\delta} x_{p}^{\delta} + \sum_{p \in P} c_{E} v_{p} \label{eq:eval_obj} \\
\text{s.t.} \quad & \sum_{\delta \in \Delta_{p}} x_{p}^{\delta} \leq \bar{y}_{p} + v_{p} && \forall p \in P \label{eq:eval_capacity} \\
& \sum_{p \in P} \sum_{\delta \in \Delta_{p}} a_{p k}^{\delta s} x_{p}^{\delta} \geq 1 && \forall k \in K \label{eq:eval_cover} \\
& v_{p} \geq 0 && \forall p \in P \\
& x_{p}^{\delta} \in \mathbb{B} && \forall p \in P, \delta \in \Delta_{p}.
\end{align}
\label{eq:evaluation}%
\end{subequations}
In essence, model (\ref{eq:evaluation}) is identical to model (\ref{eq:model}) with a single scenario, containing the actual tasks to be performed, and fixed templates. We solve this using a restricted master heuristic. First, the LP-relaxation of (\ref{eq:evaluation}) is solved to optimality using the column generation algorithm presented in Section~\ref{subsec:benders_column}. The incumbent columns are then converted back to binary variables and the model is solved using the commercial mixed-integer programming solver \texttt{CPLEX} with a time limit of 15 minutes. Computational experience suggests that this approach attains low optimality gaps.

\section{Benders Decomposition}
\label{sec:benders}

Benders decomposition is a successful variable partitioning approach for problems featuring \textit{complicating variables}: variables whose fixing results in a problem that is significantly easier to solve. In short, the Benders algorithm decomposes the problem into a Benders master problem (BMP) and one or more Benders subproblems (BSPs), the BMP containing all complicating variables and the BSPs containing all remaining variables. The BMP is solved and the values of the complicating, first-stage variables are passed to the subproblems. If a subproblem is infeasible given these first-stage decisions, a feasibility cut is added to the BMP. If the estimated objective value of the BMP underestimates the true objective value, an optimality cut is added to the BMP. The process repeats until no more cuts are generated and the original problem is solved to optimality. We refer to \cite{rahmaniani2017benders} for a recent overview on Benders decomposition.

It is clear that model (\ref{eq:model}) naturally lends itself to Benders decomposition, with template variables $\mathbf{y}$ and $\mathbf{z}$ acting as complicating variables: once template capacities have been decided, the problem decomposes into a series of independent crew scheduling problems, one per scenario. We keep the complicating variables in the BMP and construct a single BSP per scenario, containing only duty variables. In the BMP, we choose template capacities subject to all template restrictions, providing a lower bound (LB) on the optimal objective value. The BMP solution is passed to the BSP of each scenario, where it is determined how many excess duties are required given the template capacity chosen in the BMP. Since excess duties are available, this problem is always feasible and no feasibility cuts are generated. Moreover, it provides an upper bound (UB) on the optimal objective value. If the number of excess duties exceeds the estimated recovery cost of the BMP, an optimality cut is added to the BMP. This process repeats until no more cuts are identified or the gap between LB and UB is sufficiently small. As it is intractable to solve all subproblems to optimality with integral duty variables, we use a two-phase approach where we first relax duty integrality and heuristically repair the integrality later on.

The remainder of this section is structured as follows. We present the Benders master problem and Benders subproblems in Sections~\ref{subsec:benders_master} and \ref{subsec:benders_sub}, respectively. In Section~\ref{subsec:benders_column}, we present the column generation algorithm used to solve all subproblems, and in Section~\ref{subsec:benders_accelerate}, we present two acceleration techniques. In Section~\ref{subsec:benders_phases}, we describe the two-phase approach in more detail, and we conclude with a complete overview of the Benders algorithm in Section~\ref{subsec:benders_overview}. 

\subsection{Benders Master Problem}
\label{subsec:benders_master}

The Benders master problem (BMP) is a subset of model (\ref{eq:model}) and only involves decisions and constraints regarding templates. In particular, the goal is to select templates that satisfy all template restrictions and minimise a weighted sum of template costs and estimated recovery costs. An estimate of the recovery cost is provided by optimality cuts returned by the subproblems, which couple the chosen templates to the estimated recovery cost. Denote by $O_s$ the set of optimality cuts separated in scenario $s$. This set is empty in the first iteration of the algorithm. For each cut $o \in O_s$, the coefficient of template $p$ is given by $\theta_{p}^{os}$, and the constant term induced by task $k \in K_s$ is given by $\lambda_{k}^{os}$. The BMP can now be formulated as the following mixed-integer linear program:
\begin{subequations}
\begin{align}
\min \quad & \sum_{p \in P}  c_{p} y_{p} + \eta \\
\text{s.t.} \quad & y_{p} \leq C z_{p} && \forall p \in P \\
& \sum_{p \in P} z_{p} \leq \Gamma \\
& \sum_{p \in P} e_{p}^{m} y_p \leq E_m && \forall m \in M  \\ 
& \sum_{p \in P} \theta_{p}^{o s} y_p + \sum_{k \in K_s} \lambda_{k}^{o s} \leq \eta && \forall s \in S, o \in O_s \\ 
& y_{p} \in \mathbb{N} && \forall p \in P \\
& z_{p} \in \mathbb{B} && \forall p \in P.
\end{align}
\end{subequations}
The BMP is solved with a standard mixed-integer programming solver, after which the selected templates $\bar{\mathbf{y}}$ and recovery cost $\bar{\eta}$ are passed on to the Benders subproblems.

\subsection{Benders Subproblem}
\label{subsec:benders_sub}

There is a single Benders subproblem (BSP) per scenario, in which the templates $\bar{\mathbf{y}}$ chosen in the BMP are considered fixed. The goal is to determine the minimum number of excess duties needed, in addition to the capacity provided by the templates, to cover all tasks in the scenario.  Omitting scenario indices for convenience, the BSP in primal form reads:
\begin{subequations}
\begin{align}
\min \quad & \sum_{p \in P} c_{E} v_{p} \\
\text{s.t.} \quad & \sum_{\delta \in \Delta_{p}} x_{p}^{\delta} \leq \bar{y}_{p} + v_{p} && \forall p \in P \label{eq:bsp_capacity} \\
& \sum_{p \in P} \sum_{\delta \in \Delta_{p}} a_{p k}^{\delta} x_{p}^{\delta} \geq 1 && \forall k \in K \label{eq:bsp_cover} \\
& v_{p} \geq 0 && \forall p \in P \\
& x_{p}^{\delta} \geq 0 && \forall p \in P, \delta \in \Delta_{p}.
\end{align}
\label{eq:bsp_primal}%
\end{subequations}
The BSP can be readily solved using column generation (see Section~\ref{subsec:benders_column}). Due to the excess duties, the problem is always feasible and no feasibility cuts are generated. If the optimal objective value of (\ref{eq:bsp_primal}) exceeds the estimated recovery cost $\bar{\eta}$, we separate an optimality cut. To determine the expression of these optimality cuts, we first write the BSP in dual form. Introducing dual vectors $\mathbf{\theta}$ and $\mathbf{\lambda}$ for constraints (\ref{eq:bsp_capacity}) and (\ref{eq:bsp_cover}), respectively, the dual BSP reads:
\begin{subequations}
\begin{align}
\max \quad & \sum_{p \in P} \bar{y}_{p} \theta^{p} + \sum_{k \in K} \lambda_{k} \\
\text{s.t.} \quad & \theta_{p} + \sum_{k \in K} a_{p k}^{\delta} \lambda_{k} \leq 0 && \forall p \in P, \delta \in \Delta_{p} \\
& -c_E \leq \theta_{p} \leq 0 && \forall p \in P \\ 
& \lambda_{k} \geq 0 && \forall k \in K. 
\end{align}
\end{subequations}
The Benders optimality cuts can thus be written as: 
\begin{align}
	\sum_{p \in P} \theta_{p} y_{p} + \sum_{k \in K} \lambda_{k} \leq \eta.
\end{align}
Intuitively, these cuts indicate at which times of the day more template capacity is required to  reduce the recovery cost. The cuts are only valid if the BSP is solved to optimality.

\subsection{Column Generation}
\label{subsec:benders_column}

We solve all subproblems using column generation. Here, the subproblems are initialised with only a small subset of all possible duties, yielding a restricted master problem (RMP). After solving the RMP, a pricing problem is solved to identify columns with negative reduced cost that potentially improve the objective value. The algorithm repeats until no more columns with negative reduced cost exist and the subproblem is solved to optimality. The reduced cost of duty $x_{p}^{\delta}$ is given by 
\begin{align}
	-\theta_p - \sum_{k \in K} \lambda_{k} a_{pk}^{\delta}.
\end{align}
The pricing problem of identifying a duty with negative reduced cost can be decomposed over all templates. The pricing problem for a single template can be modelled as a resource-constrained shortest path problem on a suitably-defined directed acyclic graph. This graph includes all tasks that can potentially be performed within the template, as well as all arcs modelling feasible connections between consecutive tasks. A source and sink node model the departure from and arrival at the crew base, respectively. The reduced cost can be easily decomposed on the arcs in this network. Resources model the maximum duty length and maximum time without meal break. 

We start by solving all pricing problems with the heuristic labelling algorithm of \cite{breugem2022column}, augmented with a completion bound on the duty length and time without break. Once the heuristic pricing fails to return columns with negative reduced cost, we switch to the exact pricing algorithm of \cite{huisman2007column} with $\mathcal{O}(|K|^3)$ time complexity. The subproblems must be solved to optimality in order for the Benders cuts to be valid.

We apply several acceleration strategies to speed up the column generation algorithm. First, we limit the number of columns returned by any pricing problem. Second, we greedily select columns to add to the RMP based on negative reduced cost, until no more columns exist that are sufficiently task-disjoint from previously selected columns \citep{breugem2022column}. Third, we apply column management and remove columns whose reduced cost exceeds a given threshold for a certain number of iterations. Fourth, we solve all pricing problems in parallel.

\subsection{Acceleration Techniques}
\label{subsec:benders_accelerate}

In this section, we propose two acceleration techniques to speed up convergence of the Benders decomposition algorithm: valid inequalities (Section~\ref{subsubsec:benders_valid}) and the separation of Pareto-optimal cuts (Section~\ref{subsubsec:benders_pareto}).

\subsubsection{Valid Inequalities}
\label{subsubsec:benders_valid}

It is well-known that including a relaxation of the subproblem in the master problem can speed up convergence of the Benders algorithm \citep{rahmaniani2017benders}. To this end, we relax the template capacity constraints and several complicating duty rules (maximum length, meal break requirement, identical start and end depot requirement) in the crew scheduling subproblem. What remains is a min-cost flow problem, in which all tasks must be covered by flows starting from and ending at a depot. This problem can be solved in polynomial time, allowing us to efficiently compute strong lower bounds on the template capacity required at various times of the day.

Omitting scenario indices for simplicity, recall that $K$ is the set of tasks to cover. As in the pricing problem (see Section~\ref{subsec:benders_column}), we construct a directed acyclic graph containing the tasks in $K$ as well as an artificial source and sink node representing the depot. An arc in this graph models a feasible connection between two tasks.  Let $A$ be the set of arcs, and denote by $A^+(k)$ and $A^{-}(k)$ the set of outgoing and ingoing arcs of task $k \in K$, respectively. Moreover, let $A(t)$ denote the set of arcs whose traversal requires a crew member to work in the interval $[t, t+1]$. We introduce the continuous decision variable $f_a$ to represent the flow over arc $a \in A$. A lower bound on the number of required duties in the interval $[t, t+1]$ can be obtained by solving the following min-cost flow problem:
\begin{subequations}
\begin{align}
\min \quad & \sum_{a \in A(t)} f_a \label{eq:flow_obj} \\
\text{s.t.} \quad & \sum_{a \in A^{+}(k)} f_a = \sum_{a \in A^{-}(k)} f_a && \forall k \in K \label{eq:flow_balance} \\
& \sum_{a \in A^{+}(k)} f_a \geq 1 && \forall k \in K \label{eq:flow_cover} \\
& f_a \geq 0 && \forall a \in A.  \label{eq:flow_domain}
\end{align}
\label{eq:flow}%
\end{subequations}
The objective (\ref{eq:flow_obj}) is to minimise the number of arcs simultaneously used at time $t$. Constraints (\ref{eq:flow_balance}) ensure flow balance, while constraints (\ref{eq:flow_cover}) ensure that each task is covered at least once. The domain of flow variables is given by (\ref{eq:flow_domain}).

Let $V_t$ denote the optimal objective value of (\ref{eq:flow}), and let $P(t)$ be the set of templates active during the interval $[t, t+1]$. We add to the BMP the valid inequality
\begin{align}
\sum_{p \in P(t)} y_p + \frac{\eta}{c_E} \geq V_t.
\end{align}
We compute an inequality for each scenario and for values of $t$ five minutes apart, balancing computing time and the goal of appropriately covering all times of day.

\subsubsection{Pareto-Optimal Cut Selection}
\label{subsubsec:benders_pareto}

The strength of the optimality cuts has a great impact on the convergence rate of the Benders decomposition algorithm \citep{rahmaniani2017benders}. Since the crew scheduling problem is highly degenerate, however, multiple optimal dual solutions can exist which might not all be equally effective at approximating the true recovery cost function. \cite{magnanti1981accelerating} propose a method of choosing among optimal dual solutions by introducing the concept of Pareto-efficient cuts. In particular, they say an optimality cut is Pareto-efficient when, among all the solutions that are optimal at the incumbent BMP solution, it is also maximal with respect to some fixed `core point' solution to the BMP. Such a cut can be obtained by solving an auxiliary subproblem. This approach has been successfully applied to crew planning by, among others, \cite{mercier2005computational}.

Let $\bar{\gamma}$ denote the objective of the regular BSP, and let $\mathbf{w}$ denote the core point vector. The auxiliary subproblem in dual form then reads:
\begin{subequations}
\begin{align}
\max \quad & \sum_{p \in P} w_{p} \theta_{p} + \sum_{k \in K} \lambda_{k} \label{eq:pareto_obj} \\
\text{s.t.} \quad & \theta_{p} + \sum_{k \in K} a_{p k}^{\delta s} \lambda_{k} \leq 0 && \forall p \in P, \delta \in \Delta_{p} \\
&  \sum_{p \in P} \bar{y}_{p} \theta_{p} + \sum_{k \in K} \lambda_{k} = \bar{\gamma} \label{eq:pareto_constraint} \\
& -c_{E} \leq \theta_{p} \leq 0 && \forall p \in P \\ 
& \lambda_{k} \geq 0  && \forall k \in K.
\end{align}
\end{subequations}
The modified objective~(\ref{eq:pareto_obj}) is to maximise the value of the cut at the core point vector. Constraint~(\ref{eq:pareto_constraint}) ensures that the cut remains optimal at the incumbent BMP solution. Together, this ensures Pareto-efficiency of the resulting cut. 

In primal form, the changes to the dual problem result in a modification of the right-hand side of each template capacity constraint and an additional unrestricted variable $u$:
\begin{subequations}
\begin{align}
\min \quad & \sum_{p \in P} c_{E} v_{p} + \bar{\gamma} u \\
	\text{s.t.} \quad & \sum_{\delta \in \Delta_{p}} x_{p}^{\delta}  + \bar{y}^{p} u \leq w_{p} + v_{p} && \forall p \in P \\
& \sum_{p \in P} \sum_{\delta \in \Delta_{p}} a_{p k}^{\delta} x_{p}^{\delta} + u \geq 1 && \forall k \in K \\
& v_{p} \geq 0 && \forall p \in P \\
& x_{p}^{\delta} \geq 0 && \forall p \in P, \delta \in \Delta_{p}.
\end{align}
\label{eq:auxiliary_primal}%
\end{subequations}
We solve the primal auxiliary subproblem (\ref{eq:auxiliary_primal}) with the same column generation algorithm as the regular BSP. As before, the optimal dual variables of this problem are used to construct optimality cuts. Following \cite{mercier2005computational}, we initialise the core point as a constant, positive vector. In each iteration, we update the core point as $\lambda$ times the current BMP solution plus $(1 - \lambda)$ times the old core point. In line with \cite{papadakos2008practical}, we set $\lambda = 1/2$.

\subsection{Two-Phase Approach}
\label{subsec:benders_phases}

It is challenging to solve all subproblems to optimality when the duty variables are required to be integral: this would require one to incorporate the Benders decomposition algorithm into a branch-and-bound framework. To avoid this issue, we adopt a heuristic two-phase Benders decomposition approach, similar to the works of \cite{cordeau2001benders} and \cite{mercier2005computational}. In the first phase, we relax the integrality of all duty variables in the BSPs. The resulting subproblems can be readily solved to optimality using column generation. We apply Benders decomposition until some stopping criterion is met, obtaining a solution with potentially fractional duty variables in the subproblems. This solution is not necessarily feasible in the original problem and hence provides a lower bound on the optimal objective value. To repair integrality, we proceed to the second phase. Here, we iteratively fix duty variables to one in the BSPs and re-apply Benders decomposition to the resulting, restricted problem. We repeat this until all subproblem solutions are fully integral, thereby obtaining a feasible solution to the original problem and an upper bound on the optimal objective value. \cite{mercier2005computational} precede the first step by an additional step where integrality of the complicating variables in the BMP is relaxed. Since solving our BMP is not very time-consuming, we choose not to do this here. While the two-phase approach is heuristic in nature, the bounds obtained in the first and second phase allow us to compute optimality gaps. As we will see, these gaps are typically rather small. We now describe each phase in more detail.

In the first phase we relax integrality of the duty variables in the Benders subproblems, allowing us to solve the BSPs to optimality using column generation. We apply regular Benders decomposition until no more optimality cuts are generated or a time limit of one hour is reached. The first-phase solution can contain fractional duty variables, and is not necessarily a feasible solution to the original problem (\ref{eq:model}). Instead, its objective value is a lower bound on the true optimal objective value. This lower bound is generally tight, since the linear programming relaxation of the crew scheduling problem typically provides a good approximation of optimal integer solution values. 

Once the first phase terminates, we enter the second phase in which integrality of the duty variables is repaired. The second phase alternates between fixing iterations and resolving iterations. In fixing iterations, the fractional duty variable with highest value in every Benders subproblem is permanently fixed to one. After fixing, we add a valid inequality to the BMP stating that the selected templates provide capacity for at least all fixed duties. In resolving iterations, several iterations of Benders decomposition are applied to allow the BMP solution to change in response to the fixing: it might happen that additional templates must be selected in response to duty variables being rounded to one. We repeat the sequence of fixing and resolving iterations until all duty variables are integral and the solution is feasible in model (\ref{eq:model}). The objective value of this solution provides an upper bound on the optimal objective value. The Benders lower bound obtained in the second phase is no longer valid, as fixing decisions are not necessarily optimal. 

\subsection{Overview}
\label{subsec:benders_overview}

Algorithm~\ref{alg:benders} provides a brief overview of the complete two-phase accelerated Benders decomposition algorithm. We decide to solve the BMP heuristically by imposing a time limit of ten seconds to the mixed-integer programming solver. Whenever no more cuts are generated and optimality is required, we iteratively increase the BMP time limit to obtain provably optimal solutions.

\begin{algorithm}[htbp!]
\caption{Two-phase accelerated Benders decomposition}
\label{alg:benders}
\begin{algorithmic}
\State Phase $\gets$ I
\State Add valid inequalities to BMP \Comment{Section~\ref{subsec:benders_accelerate}}
\State Initialise core point vector $\mathbf{w}$
\While{Phase $=$ I or solution not integral}
	\State Solve the BMP to obtain $(\bar{\mathbf{y}}, \bar{\eta})$  \Comment{Section~\ref{subsec:benders_master}}
	\For{Each scenario} 
		\State Solve a BSP given $\bar{\mathbf{y}}$ \Comment{Sections~\ref{subsec:benders_sub}-\ref{subsec:benders_column}}
		\If{Objective exceeds $\bar{\eta}$}
			\State Solve auxiliary BSP given $\mathbf{w}$ \Comment{Section~\ref{subsec:benders_accelerate}}
			\State Add optimality cut to BMP
		\EndIf
	\EndFor
	\If{Phase $=$ I, and no cuts or time limit reached}
		\State Phase $\gets$ II
	\EndIf
	\If{Phase $=$ II, and no cuts or iteration limit reached}
		\For{Each scenario}
			\State Fix one or more duties in BSP \Comment{Section~\ref{subsec:benders_phases}}
		\EndFor
	\EndIf
	\State $\mathbf{w} \gets \lambda \bar{\mathbf{y}} + (1 - \lambda) \mathbf{w}$
\EndWhile
\end{algorithmic}
\end{algorithm}

\section{Benchmark}
\label{sec:benchmark}

We benchmark our Benders decomposition algorithm against the column generation heuristic of \cite{rahlmann2021robust}, who studied a highly similar robust tactical crew scheduling problem. Their problem setting differs from ours in two important ways. First, they do not assume a pre-specified list of template types, but allow templates to start at all possible times of the day. Second, they do not impose rostering restrictions on the set of chosen templates, i.e., they do not impose constraints (\ref{eq:template_activation})-(\ref{eq:template_mix}).

We provide a brief overview of the benchmark algorithm in Section~\ref{subsec:benchmark_overview}, and refer to \cite{rahlmann2021robust} for a detailed description. We describe several changes in implementation in Section~\ref{subsec:benchmark_changes}, and compare the performance of the two algorithms on instances with pre-specified templates and no template restrictions in Section~\ref{subsec:computational_benchmark}. 

\subsection{Overview}
\label{subsec:benchmark_overview}

Whereas we assume that templates must be chosen from a pre-specified list, \cite{rahlmann2021robust} consider a slightly more flexible setting where templates are allowed to start at arbitrary times of the day. In particular, a template is defined as a compatible sequence of duties across multiple scenarios, duties being compatible when their start and end times comply with the maximum template length restriction. In other words, each template specifies a duty to be performed in each scenario. Since this leads to a huge number of possible templates, a column generation algorithm is proposed to solve the problem. 

This algorithm alternates between two different pricing problems. The first pricing problem consists of aggregating duties to feasible templates with negative reduced cost. To this end, a tailored pricing network is proposed. It contains a partition for each scenario, a node for each duty in a scenario, and an artificial source and sink node. Duties in different layers are connected when they can be performed in a single template. Each feasible template corresponds to a source-sink path in this network, and a labelling algorithm is used to identify negative reduced cost templates. Since the number of possible duties is large, the network is initialised with a small subset of duties only. In particular, an ordinary crew scheduling problem is solved for each scenario, and the solution to this problem is used as initial node set. 

When the first pricing problem fails to return new templates, a second pricing problem is invoked to identify promising duties whose addition to the network might allow for the generation of new negative reduced cost templates. This pricing problem decomposes per scenario, and is highly similar in structure to the pricing problem described in Section~\ref{subsec:benders_column}. However, the reduced cost criterion accounts for the fact that duties must be aggregated to templates later on. Finally, a fixing rule is used to obtain integer solutions.

\subsection{Adaptations in Implementation}
\label{subsec:benchmark_changes}

We make several modelling choices to ensure that the solution space of the benchmark algorithm is the same as that of our Benders decomposition algorithm. In contrast with \cite{rahlmann2021robust}, we assume a single crew base and single day of operation, and omit capacity constraints from consideration (though these can easily be incorporated in our approach). Since our problem setting features a pre-specified list of allowed template starting times, we artificially modify the source arcs in the tailored network to align with these starting times. We use the same parameter settings as in the benchmark paper, and use the same pricing algorithms as for the Benders decomposition algorithm. 

The number of possible duties in our passenger railway application is considerably larger than in freight railway, leading to a rapid growth of the tailored network throughout the course of the column generation algorithm. In response, we have made several minor changes in implementation to limit the size of the network. First, we only connect the source node to duties in the first scenario. Without loss of optimality, this reduces the number of source arcs and imposes that all templates use exactly one duty in each scenario. Second, we apply a task-disjointness filter (see Section~\ref{subsec:benders_column}) before adding duties from the second pricing problem to the network. Third, we remove duties based on their reduced cost as soon as the network contains more than a specified number of nodes. In our experiments, we use a maximum size of 5,000.  

\section{Computational Performance}
\label{sec:computational}

In this section, we evaluate the computational performance of the proposed Benders algorithm by conducting extensive experiments on real-life instances from NS. We describe the instances and parameter settings in Sections~\ref{subsec:computational_instances} and \ref{subsec:computational_parameter}, respectively.  In Section~\ref{subsec:computational_acceleration}, we analyse the effect of the various acceleration techniques of Section~\ref{subsec:benders_accelerate} on the first phase performance, and in Section~\ref{subsec:computational_phases}, we consider the overall performance of the two-phase algorithm. We compare our algorithm with that of \cite{rahlmann2021robust} in Section~\ref{subsec:computational_benchmark}. 

\subsection{Instances}
\label{subsec:computational_instances}

We use seven weeks of historical planning data of NS, spanning the period of October 4 to November 21 of 2021. For each day, the detailed data contains the set of tasks to be performed and the scheduled duties for all train guards in the Netherlands. The first three weeks are used to construct scenario sets, i.e., as input in the tactical planning phase, while the last four weeks are used for evaluation of the operational planning phase. 

\begin{table}[htbp!]
  	\centering 
  	\renewcommand{\arraystretch}{1}
  	\caption{Instance characteristics. For each crew base and day of week, the number of tasks and similarity ratio are given. Results are averaged over all weeks.}
  	\label{table:instances} 
  	\begin{adjustbox}{max width= \textwidth}
  	\begin{tabular}{@{\extracolsep{0pt}}lrrrrrrrrrrrrrr} 
		\toprule 
		Base & \multicolumn{2}{c}{Mon} & \multicolumn{2}{c}{Tue} & \multicolumn{2}{c}{Wed} & \multicolumn{2}{c}{Thu} & \multicolumn{2}{c}{Fri} & \multicolumn{2}{c}{Sat} & \multicolumn{2}{c}{Sun} \\
  		\midrule
		Ah  & 507 & 0.69 & 500 & 0.64 & 471 & 0.83 & 576 & 0.73 & 552 & 0.72 & 460 & 0.51 & 349 & 0.34 \\
Amf & 499 & 0.61 & 499 & 0.59 & 465 & 0.69 & 509 & 0.66 & 496 & 0.51 & 366 & 0.13 & 345 & 0.16 \\
Asd & 859 & 0.71 & 839 & 0.69 & 893 & 0.79 & 940 & 0.74 & 910 & 0.62 & 669 & 0.21 & 560 & 0.17 \\
Ehv & 593 & 0.60 & 612 & 0.58 & 576 & 0.68 & 641 & 0.63 & 587 & 0.55 & 502 & 0.36 & 418 & 0.31 \\
Gvc & 898 & 0.79 & 917 & 0.80 & 909 & 0.86 & 948 & 0.77 & 914 & 0.71 & 762 & 0.42 & 631 & 0.33 \\
Nm  & 430 & 0.68 & 440 & 0.70 & 429 & 0.84 & 466 & 0.76 & 496 & 0.74 & 382 & 0.37 & 307 & 0.37 \\
Rtd & 615 & 0.74 & 598 & 0.79 & 552 & 0.81 & 655 & 0.81 & 633 & 0.70 & 502 & 0.41 & 386 & 0.25 \\
Ut  & 843 & 0.73 & 814 & 0.73 & 801 & 0.77 & 838 & 0.76 & 891 & 0.65 & 744 & 0.35 & 581 & 0.25 \\
\midrule
\textbf{Avg.} & \textbf{656} & \textbf{0.69} & \textbf{652} & \textbf{0.69} & \textbf{637} & \textbf{0.79} & \textbf{697} & \textbf{0.73} & \textbf{685} & \textbf{0.65} & \textbf{549} & \textbf{0.35} & \textbf{447} & \textbf{0.27} \\ 
\bottomrule 
\end{tabular} 
\end{adjustbox}
\end{table}

Using the historic assignment of tasks to crew bases through duties, we construct $7 \times 8$ instances, one for each combination of weekday and one of eight crew bases. Table~\ref{table:instances} lists, for each instance, the average number of tasks and the similarity ratio. The similarity ratio, defined as the average fraction of tasks that appears in two consecutive planning weeks, is a rough proxy of the variability in task sets. The instances range in size: Crew bases Amersfoort (Amf) and Nijmegen (Nm) contain between 300 and 500 tasks per day, while crew bases The Hague (Gvc) and Amsterdam (Asd) can contain well over 900 tasks per day. Generally, more tasks are performed throughout the week than in the weekend. The task sets appear to be more stable throughout the week than in the weekend: The average similarity ratio on Wednesday equals 0.79, but drops to 0.27 on Sunday. This can be explained by the fact that most maintenance works and events are scheduled in the weekend, necessitating many changes to the timetable.

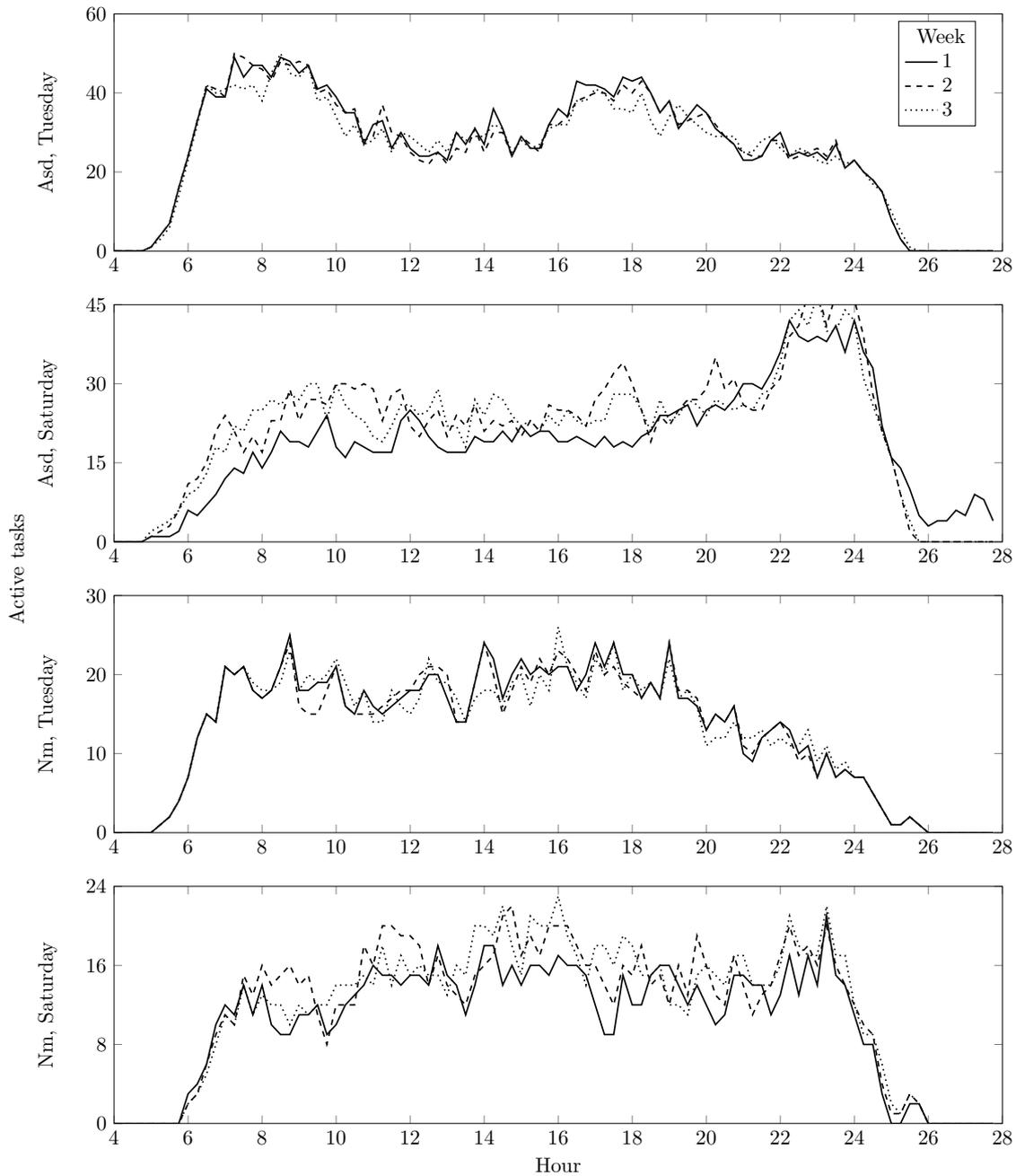
\begin{figure}[htbp!]
\centering
 \begin{adjustbox}{max width= \textwidth}
\begin{tikzpicture}
\begin{groupplot}[group style={group name=my plots,group size= 1 by 4,horizontal sep =1.5cm, vertical sep =1cm}, width=18cm, height=6cm, groupplot xlabel={Hour}, groupplot ylabel={Active tasks}]
       \nextgroupplot[
	ylabel={Asd, Tuesday},
         ymin =0,     
	 ymax=60,  
	 xmin=4, 
	xmax=28,
         minor x tick num = 0,
         minor tick length=2ex,
         ytick = {0, 20, 40, 60},
	legend pos = north east,
	legend cell align = left,
        ]
        \addlegendimage{empty legend}
        \addplot +[style={black, mark=none, thick}] 
table [col sep= semicolon, x index = 0, y index=1] {\subfix{taskScenarios_Asd_Tue.csv}};

 	\addplot +[style={black, dashed, mark=none, thick}]
table [col sep= semicolon, x index = 0, y index=2] {\subfix{taskScenarios_Asd_Tue.csv}};

	\addplot +[style={black, dotted, mark=none, thick}]
table [col sep= semicolon, x index = 0, y index=3] {\subfix{taskScenarios_Asd_Tue.csv}};

	\addlegendentry{\hspace{-.6cm} Week}
   	\addlegendentry{1}
   	\addlegendentry{2}
	\addlegendentry{3}

\nextgroupplot[
         ymin =0,     
	   ymax=45,  
	   xmin=4,
	   xmax=28,
	ylabel={Asd, Saturday},
        minor x tick num = 0,
        minor tick length=2ex,
         ytick = {0, 15, 30, 45},
	legend pos = north west,
	legend cell align = left,
        ]
        \addlegendimage{empty legend}
        \addplot +[style={black, mark=none, thick}] 
table [col sep= semicolon, x index = 0, y index=1] {\subfix{taskScenarios_Asd_Sat.csv}};

 	\addplot +[style={black, dashed, mark=none, thick}]
table [col sep= semicolon, x index = 0, y index=2] {\subfix{taskScenarios_Asd_Sat.csv}};

	\addplot +[style={black, dotted, mark=none, thick}]
table [col sep= semicolon, x index = 0, y index=3] {\subfix{taskScenarios_Asd_Sat.csv}};

\nextgroupplot[
         ymin =0,     
	   ymax=30,  
	   xmin=4,
	   xmax=28,
	ylabel={Nm, Tuesday},
        minor x tick num = 0,
        minor tick length=2ex,
         ytick = {0, 10, 20, 30},
	legend pos = north west,
	legend cell align = left,
        ]
        \addlegendimage{empty legend}
        \addplot +[style={black, mark=none, thick}] 
table [col sep= semicolon, x index = 0, y index=1] {\subfix{taskScenarios_Nm_Tue.csv}};

 	\addplot +[style={black, dashed, mark=none, thick}]
table [col sep= semicolon, x index = 0, y index=2] {\subfix{taskScenarios_Nm_Tue.csv}};

	\addplot +[style={black, dotted, mark=none, thick}]
table [col sep= semicolon, x index = 0, y index=3] {\subfix{taskScenarios_Nm_Tue.csv}};

\nextgroupplot[
         ymin =0,     
	   ymax=24,  
	   xmin=4,
	   xmax=28,
	ylabel={Nm, Saturday},
        minor x tick num = 0,
        minor tick length=2ex,
         ytick = {0, 8, 16, 24},
	legend pos = north west,
	legend cell align = left,
        ]
        \addlegendimage{empty legend}
        \addplot +[style={black, mark=none, thick}] 
table [col sep= semicolon, x index = 0, y index=1] {\subfix{taskScenarios_Nm_Sat.csv}};

 	\addplot +[style={black, dashed, mark=none, thick}]
table [col sep= semicolon, x index = 0, y index=2] {\subfix{taskScenarios_Nm_Sat.csv}};

	\addplot +[style={black, dotted, mark=none, thick}]
table [col sep= semicolon, x index = 0, y index=3] {\subfix{taskScenarios_Nm_Sat.csv}};
\end{groupplot}
\end{tikzpicture}
\end{adjustbox}
\caption{Temporal distribution of tasks for crew bases Amsterdam (Asd) and Nijmegen (Nm) on Tuesdays and Saturdays in the first three weeks of planning data.}
\label{fig:scenarios}
\end{figure}

To gain more insight into the scenario sets, Figure~\ref{fig:scenarios} presents the number of active tasks throughout the day for three scenarios of crew bases Amsterdam (Asd) and Nijmegen (Nm) on Tuesday and Sunday. The number of active tasks is computed as the number of tasks whose start and end time are before and after a given time, respectively. We observe a strong heterogeneity in capacity demand patterns across crew bases and days of the week. Nonetheless, the pattern of capacity demand is generally quite similar across weeks. This suggests that robust template capacities might exist, i.e., templates that are able to accurately provide capacity in all weeks of the planning period. The relatively low similarity ratios on Saturday (see Table~\ref{table:instances}) are reflected in Figure~\ref{fig:scenarios}, showing more variety on the Saturdays than on the Tuesdays. As a result, providing robust schedules for the weekend is expected to be more challenging.

\subsection{Parameter Settings}
\label{subsec:computational_parameter}

Unless indicated otherwise, we use the following parameter settings. All regular templates span 9.5 hours and are allowed to start every 30 minutes. There is a single reserve template that spans the whole day. The maximum duty length equals nine hours. We consider problems with three scenarios, in which at most 15 unique template types and 10\% reserve templates may be selected. To reflect the different demand patterns, the maximum percentage of early and late templates depends on the day of the week: from Monday to Wednesday we allow at most 30\% early and 30\% late templates, on Thursday and Friday 30\% and 45\%, and on Saturday and Sunday at most 20\% and 55\%.

A single template costs 10,000, whereas an excess duty costs 40,000. In the evaluation phase, we use a variable cost of one unit per second, i.e., each worked hour contributes 3,600 towards the objective. 

We use the following column generation parameters. At most 50 columns are returned per pricing problem, and only columns that are at least 50\% task-disjoint from previously selected columns are added to the RMP. We remove columns whose reduced cost was positive in the last 15 iterations, and solve all pricing problems in parallel on eight threads. All experiments are conducted on a computing server with 16GB RAM and an AMD Rome 7H12 processor. CPLEX 22.1.0 is used to solve all (mixed-integer) linear programs.

\subsection{Effect of Acceleration Techniques}
\label{subsec:computational_acceleration}

To analyse the effect of acceleration techniques on the performance of the Benders decomposition algorithm, we conduct experiments on the Monday instances where we only consider the first phase of the algorithm and impose a time limit of three hours. We compare three versions of the Benders algorithm: a default version, a version with Pareto-efficient cuts, and a version with Pareto-efficient cuts and valid inequalities. Table~\ref{table:effect_cuts} reports the results of these experiments. For each parameter setting and version of the algorithm, it reports the average number of iterations, computing time in seconds, split up over time spent on the BMP and time spent on the BSPs, the optimality gap in the final iteration, and the number of instances solved to optimality. Since we consider only the first phase of the algorithm, where duty variables are allowed to be fractional, optimality in this experiment does not necessarily imply global optimality in the original problem.

\begin{table}[htbp!]
  	\centering 
  	\renewcommand{\arraystretch}{1}
  	\caption{Effect of acceleration techniques on computational performance of first phase Benders decomposition for Monday instances.}
  	\label{table:effect_cuts} 
  	\begin{adjustbox}{max width= \textwidth}
  	\begin{tabular}{@{\extracolsep{0pt}}lllrrrrrr} 
		\toprule 
		& & & & \multicolumn{3}{c}{Time (s)} \\
		\cmidrule(l){5-7} 
  		Method & Types & Scenarios & Iterations & BMP & BSP & Total & Gap (\%) & Optimal \\
  		\midrule
		\multirow{10}{*}{Default} & \multirow{3}{*}{10} & 1 & 458.5 & 4,758.5 & 1,531.6 & 6,290.1 & 0.15 & 7/8 \\
		 & & 2 & 351.9 & 3,388.4 & 1,862.4 & 5,250.8 & 1.34 & 6/8 \\
		 & & 3 & 529.4 & 3,904.3 & 4,443.4 & 8,437.6 & 1.27 & 4/8 \\
		\cmidrule{2-9}
		 & \multirow{3}{*}{15} & 1 & 596.0 & 2,787.8 & 2,292.6 & 5,080.4 & 0.15 & 7/8 \\
		 & & 2 & 584.0 & 1,267.4 & 4,294.9 & 5,562.3 & 1.14 & 6/8 \\
		 & & 3 & 593.3 & 1,373.5 & 5,282.1 & 6,655.6 & 1.07 & 5/8 \\
		\cmidrule{2-9}
		 & \multirow{3}{*}{20} & 1 & 621.4 & 2,234.4 & 2,140.3 & 4,374.6 & 0.50 & 7/8 \\
		 & & 2 & 703.0 & 1,306.5 & 4,774.0 & 6,080.5 & 0.22 & 7/8 \\
		 & & 3 & 574.3 & 1,335.1 & 5,460.6 & 6,795.8 & 2.33 & 4/8 \\
		\cmidrule{2-9} 
		& \textbf{Avg.} & & \textbf{556.8} & \textbf{2,484.0} & \textbf{3,564.7} & \textbf{6,048.6} & \textbf{0.91} & \textbf{53/72} \\ 
		\midrule
		\multirow{10}{*}{Pareto} & \multirow{3}{*}{10} & 1 & 192.8 & 2,851.9 & 1,309.9 & 4,161.8 & 0.15 & 7/8 \\
		 & & 2 & 186.0 & 2,860.3 & 2,384.1 & 5,244.4 & 0.15 & 7/8 \\
		 & & 3 & 130.6 & 1,885.0 & 1,761.1 & 3,646.1 & 0.08 & 7/8 \\
		\cmidrule{2-9}
		 & \multirow{3}{*}{15} & 1 & 160.4 & 1,513.8 & 983.4 & 2,497.1 & 0 & 8/8 \\
		 & & 2 & 126.4 & 1,012.9 & 1,555.0 & 2,567.9 & 0.02 & 7/8 \\
		 & & 3 & 85.8 & 507.5 & 1,313.5 & 1,821.0 & 0 & 8/8 \\
		\cmidrule{2-9}
		 & \multirow{3}{*}{20} & 1 & 188.3 & 1,290.5 & 969.6 & 2,260.1 & 0 & 8/8 \\
		 & & 2 & 93,1 & 869.0 & 951.9 & 1,820.9 & 0 & 8/8 \\
		 & & 3 & 88.5 & 619.4 & 1,472.4 & 2,091.8 & 0 & 8/8 \\
		\cmidrule{2-9} 
		& \textbf{Avg.} & & \textbf{139.1} & \textbf{1,490.0} & \textbf{1,411.2} & \textbf{2,901.2} & \textbf{0.04} & \textbf{68/72} \\ 
		\midrule
		\multirow{10}{*}{Pareto + VI} & \multirow{3}{*}{10} & 1 & 183.1 & 3,424.3 & 1,090.3 & 4,514.5 & 0.20 & 6/8 \\
		 & & 2 & 147.8 & 3,109.4 & 1,736.9 & 4,846.3 & 0.08 & 7/8 \\
		 & & 3 & 115.5 & 2,296.4 & 1,811.4 & 4,107.8 & 0.08 & 7/8 \\
		\cmidrule{2-9}
		 & \multirow{3}{*}{15} & 1 & 105.8 & 712.5 & 777.8 & 1,490.3 & 0 & 8/8 \\
		 & & 2 & 105.0 & 963.8 & 1,379.5 & 2,343.3 & 0 & 8/8 \\
		 & & 3 & 73.8 & 523.8 & 1.052.8 & 1,576.5 & 0 & 8/8 \\
		\cmidrule{2-9}
		 & \multirow{3}{*}{20} & 1 & 135.9 & 1,268.6 & 840.6 & 2,109.3 & 0 & 8/8 \\
		 & & 2 & 102.0 & 704.9 & 1,219.1 & 1,924.0 & 0 & 8/8 \\
		 & & 3 & 70.6 & 427.9 & 986.1 & 1,414.0 & 0 & 8/8 \\
		\cmidrule{2-9} 
		& \textbf{Avg.} & & \textbf{115.5} & \textbf{1,492.4} & \textbf{1,210.5} & \textbf{2,702.9} & \textbf{0.04} & \textbf{68/72} \\ 
		\bottomrule 
	\end{tabular} 
  	\end{adjustbox}
\end{table}

Table~\ref{table:effect_cuts} shows that separating Pareto-efficient cuts has a drastic positive impact on the computational performance. Although the time per iteration goes up slightly, the algorithm requires only one fourth of the number of iterations, less than half of the computing time, and attains almost a percentage point reduction in final optimality gap. Adding valid inequalities leads to a reduction in average computing time of about seven percent. This effect is mostly driven by a sharp reduction in the number of 
iterations, as the time spent on the BMP per iteration expectedly goes up due to the extra constraints. The valid inequalities mainly lead to reduced computing times on instances that were already solved to optimality, as the final optimality gap and number of instances solved to optimality does not increase compared to the version with only Pareto-efficient cuts. We observe that non-optimality is often caused by a small optimality gap in the BMP that could not be closed by the solver within the time limit. Nonetheless, we observe consistent speed-ups across all instances.

\begin{figure}[htbp!]
\centering
\begin{tikzpicture}
\begin{axis}[height=8cm,
		 xlabel={Time (min)},
		ylabel={Fraction of best known LB},
             width=14cm,
        	ymin=0.75, 
		ymax=1.75, 
		xmin=0,
		xmax=60, 
		extra x ticks={0},
	         grid,	
	         grid style={opacity = 0.5, dotted, black}, 
		disabledatascaling, 
		legend cell align={left},
		xtick={15, 30, 45, 60}, 
		ytick={0.75,  1.00, 1.25, 1.5, 1.75, 2.0}]
        
    \addplot[very thick, dashed, color=black!30] table[x index=0, y index=1, col sep=semicolon] {\subfix{progressionBounds.csv}};

    \addplot[very thick, dashed, color=black]
     table[x index=0, y index=2, col sep=semicolon] {\subfix{progressionBounds.csv}};
   
    \addplot[very thick, color=black!30] table[x index=0, y index=3, col sep=semicolon]{\subfix{progressionBounds.csv}};

    \addplot[very thick, color=black] table[x index=0, y index=4, col sep=semicolon] {\subfix{progressionBounds.csv}};

    \addplot[very thick, dotted, color=black!30] table[x index=0, y index=5, col sep=semicolon]{\subfix{progressionBounds.csv}};

    \addplot[very thick, dotted, color=black] table[x index=0, y index=6, col sep=semicolon] {\subfix{progressionBounds.csv}};

\legend{Default (LB), Default (UB), Pareto (LB), Pareto (UB), Pareto + VI (LB), Pareto + VI (UB)};
\end{axis}
\end{tikzpicture}
\caption{Effect of Pareto-optimal cuts and valid inequalities on progression of first phase lower and upper bounds.}
\label{fig:results_bounds}
\end{figure}
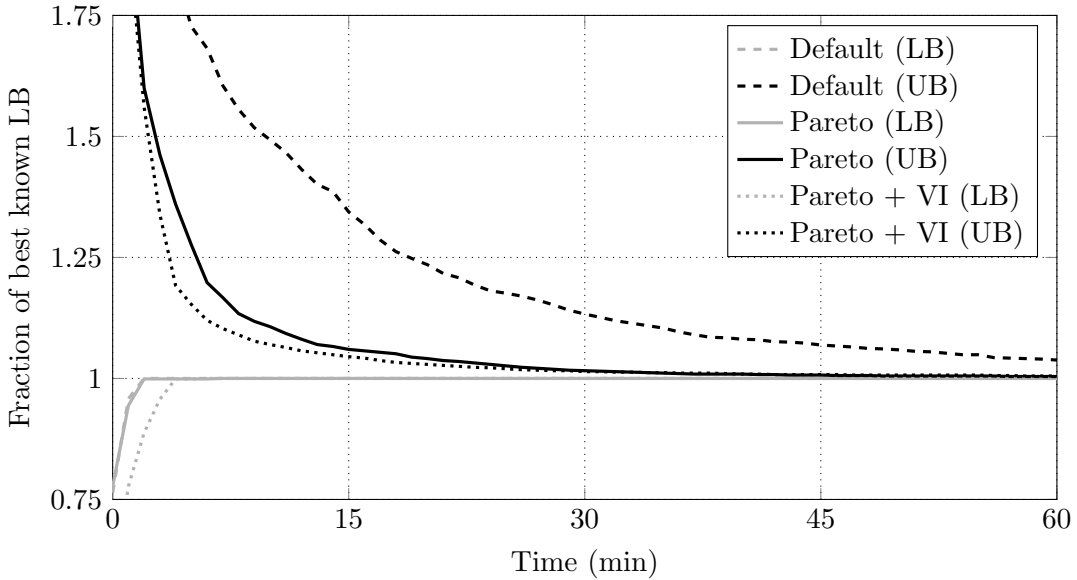

To better understand the impact of the acceleration techniques, Figure~\ref{fig:results_bounds} displays the development of the first phase Benders upper and lower bounds of the three versions of the algorithm. It shows the LB and UB as fraction of the best known LB, averaged over all instances and truncated to the first hour to highlight differences between the various approaches. It is clear that the performance improvement is purely driven by better UB behaviour. While the lower bound is often already at its optimal value after a single iteration in all approaches, the acceleration techniques lead to much faster convergence of the upper bound. The Pareto-efficient cuts have a particularly large impact, as they ensure that the recovery costs are better estimated through the optimality cuts. The valid inequalities lead to a small further improvement, by guaranteeing that solutions are close to those that require few or no excess duties. The lower bound of the version with valid inequalities slightly lags behind the other two versions, as computing the valid inequalities takes a small amount of time.

\subsection{Performance Two-Phase Approach}
\label{subsec:computational_phases}

We now analyse the performance of the complete two-phase algorithm with both acceleration techniques. We impose a time limit of one hour for the first phase, and consider all instances in our sample. Table~\ref{table:phases} reports the results. For all instances and each of the two phases, it reports the average number of iterations, computing time, and optimality gap. The first-phase gap is computed using the first-phase lower and upper bound, whereas the second-phase gap is computed using the first-phase lower bound and second-phase upper bound (see Section~\ref{subsec:benders_phases} for more detail). In addition, we report the maximum second-phase optimality gap over all instances and number of instances solved to optimality.

\begin{table}[htbp!]
  	\centering 
  	\renewcommand{\arraystretch}{1}
  	\caption{Computational performance of two-phase accelerated Benders decomposition.}
  	\label{table:phases} 
  	\begin{adjustbox}{max width= \textwidth}
  	\begin{tabular}{@{\extracolsep{0pt}}llrrrrrrrr} 
		\toprule 
		& & \multicolumn{4}{c}{First phase} & \multicolumn{3}{c}{Second phase} \\ 
		\cmidrule(l){3-5} \cmidrule(l){6-9} 
  		Types & Scenarios & Iterations & Time (s) & Gap (\%) & Iterations & Time (s) & Gap (\%) & Max gap (\%) & Optimal \\
  		\midrule
		\multirow{3}{*}{10} & 1 & 143.5 & 1,966.4 & 0.28 & 19.3 & 165.5 & 0.82 & 8.33 & 38/56 \\
		 & 2 & 113.5 & 1,933.9 & 0.43 & 34.2 & 233.5 & 1.78 & 43.75 & 31/56 \\
		 & 3 & 109.0 & 2,170.1 & 0.65 & 34.6 & 166.4 & 1.21 & 6.82 & 29/56 \\
		\midrule
		 \multirow{3}{*}{15} & 1 & 131.5 & 1,283.9 & 0.14 & 23.7 & 139.5 & 0.36 & 3.24 & 43/56 \\
		  & 2 & 104.1 & 1,357.6 & 0.08 & 35.3 & 151.8 & 0.37 & 4.00 & 43/56 \\
		  & 3 & 88.5 & 1,476.1 & 0.22 & 44.8 & 259.9 & 0.81 & 8.26 & 34/56 \\
		\midrule
		  \multirow{3}{*}{20} & 1 & 131.5 & 1,325.6 & 0.06 & 22.3 & 97.7 & 0.39 & 3.08 & 43/56 \\
		  & 2 & 95.2 & 1,290.5 & 0.14 & 34.3 & 123.1 & 0.58 & 5.26 & 40/56 \\
		  & 3 & 78.4 & 1,262.6 & 0.13 & 43.1 & 168.9 & 0.96 & 3.92 & 29/56 \\
		\midrule
		\textbf{Avg.} & & \textbf{110.6} & \textbf{1,563.0} & \textbf{0.24} & \textbf{32.4} & \textbf{167.3} & \textbf{0.81} & \textbf{11.99} & \textbf{330/504} \\ 	 
		\bottomrule 
	\end{tabular} 
  	\end{adjustbox}
\end{table}

We find that most computing time is spent on the first phase. The number of iterations and time per iteration is higher there, indicating that the variable fixing in the second phase restricts the solution space and reduces computing time. A large majority of instances is solved to optimality in the first phase, where a very low average optimality gap of 0.24\% is attained. The optimality gap increases to an average 0.81\% due to the heuristic fixing decisions in the second phase. Still, over 65\% of instances is solved to global optimality, highlighting the success of our two-phase approach. This is further confirmed by the fact that all second-phase optimality gaps are well below ten percent, except for a single outlier at 43.75\% due to a sequence of ill-performing fixing decisions. Interestingly, the number of iterations decreases as the number of scenarios grows. This might arise from the fact that more cuts are separated per iteration, speeding up convergence of the Benders algorithm. At the same time, the computing time per iteration and optimality gap slightly increase as the number of scenarios grows. Nonetheless, average computing times below half an hour are more than reasonable for a tactical planning problem.

\subsection{Comparison with Benchmark}
\label{subsec:computational_benchmark}

We benchmark our Benders decomposition algorithm against the column generation algorithm of \cite{rahlmann2021robust}. To this extent, we alter our instances and remove all template rostering constraints. In particular, we impose no restriction on the maximum number of template types or maximum number of early and late templates. We forbid the use of reserve templates, and enforce that all templates start at the half hour in both methods. In case the benchmark algorithm outputs templates with length below 9.5 hours, we articially expand these before entering the operational planning phase. 

\begin{table}[htbp!]
  	\centering 
  	\renewcommand{\arraystretch}{1}
  	\caption{Results Benders decomposition and benchmark algorithm for the first-stage problem.}
  	\label{table:benchmark_tactical} 
  	\begin{adjustbox}{max width= \textwidth}
  	\begin{tabular}{@{\extracolsep{0pt}}lrrrrrr} 
		\toprule 
		& \multicolumn{3}{c}{Benders} & \multicolumn{3}{c}{Benchmark} \\
		\cmidrule(l){2-4} \cmidrule(l){5-7} 
  		Scenarios & Time (s) & Gap (\%) & Optimal & Time (s) & Gap (\%) & Optimal \\
  		\midrule
		1 & 1,252.7 & 0.39 & 49/56 & 48.9 & 0.72 & 30/56 \\
		2 & 1,282.2 & 0.55 & 42/56 & 252.7 & 1.04 & 25/56 \\
		3 & 1,185.4 & 0.62 & 34/56 & 582.8 & 1.62 & 11/56 \\ 
		\bottomrule 
	\end{tabular} 
  	\end{adjustbox}
\end{table}

Table~\ref{table:benchmark_tactical} reports the computational performance of both methods. For each number of scenarios and method, it shows the average computing time, optimality gap, and number of instances solved to optimality. While the benchmark is on average significantly faster than our method, the Benders decomposition algorithm performs best in terms of quality for all sizes of the scenario set. The difference is especially large for the case of interest with three scenarios, where our method attains an optimality gap that is one full percentage point lower and solves 23 more instances to optimality. The two algorithms display different scaling behaviour, both in terms of time and quality. While the Benders algoritm shows stable performance with respect to the number of scenarios, the running times and optimality gaps of the column generation algorithm grow rapidly as the scenario set grows.

\begin{table}[htbp!]
  	\centering 
  	\renewcommand{\arraystretch}{1}
  	\caption{Results Benders decomposition and benchmark algorithm for the second-stage problem.}
  	\label{table:benchmark_second} 
  	\begin{adjustbox}{max width= \textwidth}
  	\begin{tabular}{@{\extracolsep{0pt}}lrrrrrrrr} 
		\toprule 
		& \multicolumn{4}{c}{Benders} & \multicolumn{4}{c}{Benchmark} \\
		\cmidrule(l){2-5} \cmidrule(l){6-9} 
  		Scenarios & Templates & Duties & Excess & Avg. workload (h) & Templates & Duties & Excess & Avg. workload (h) \\
  		\midrule
		1 & 60.07 & 62.47 & 3.03 & 7.68 & 60.27 & 61.81 & 2.25 & 7.76 \\
		2 & 61.63 & 62.40 & 1.51 & 7.70 & 61.84 & 61.97 & 1.08 & 7.73 \\
		3 & 62.09 & 62.48 & 1.25 & 7.68 & 62.86 & 62.41 & 0.80 & 7.67 \\ 
		\bottomrule 
	\end{tabular} 
  	\end{adjustbox}
\end{table}

Table~\ref{table:benchmark_second} shows the solution quality of both approaches as evaluated in the operational planning phase. For each combination of solution approach and number of scenarios, it lists the number of selected templates, number of duties and excess duties in the operational phase, and average workload in hours. All results are averaged over all instances and weeks of the evaluation period. We observe that the improved solution quality in the tactical planning phase is somewhat penalised in the operational planning phase. Since the Benders algorithm selects fewer templates in the first stage, it requires on average more duties and excess duties in the second stage. For three scenarios, however, the total number of crew members required, obtained by summing template capacity and excess duties, is approximately equal for both methods. As we will also see in Section~\ref{sec:practical}, the results in Table~\ref{table:benchmark_second} mainly show that a finite scenario set will never yield solutions that are fully robust, i.e., require no excess duties. 

\section{Practical Insights}
\label{sec:practical}

In this section, we provide several practical insights by analysing the effect of various restrictions on solution quality. Throughout the following, we will evaluate the quality of chosen templates by analysing the number of templates chosen in the tactical phase, the number of regular and excess duties required in the operational phase, and the average workload of the operational duties. Unless indicated otherwise, all results are averaged over the eight crew bases, seven days of the week, and four weeks used for evaluation purposes. We use the same instances and parameters as in Section~\ref{sec:computational}, and choose templates using the Benders decomposition algorithm with both acceleration techniques. 

We analyse the effect of the number of scenarios and template types in Section~\ref{subsec:practical_scenarios}, and present results on the number of reserves in Section~\ref{subsec:practical_reserves}. In Section~\ref{subsec:practical_length}, we present the effect of the template length on solution quality, and we split up results by day of the week in Section~\ref{subsec:practical_day}.

\subsection{Number of Scenarios and Types}
\label{subsec:practical_scenarios}

Table~\ref{table:effect_scenarios} displays the effect of the number of scenarios and the maximum number of template types on the solution quality. Increasing the number of scenarios generally leads to more robust solutions. With a larger scenario set, more templates are chosen to cover the tasks in all scenarios, leading to a reduction in the number of excess duties required in the operational phase. For example, for at most 15 unique template types, moving from one to three scenarios means that, on average, 1.51 additional templates are selected and 1.60 fewer excess duties are required. This is beneficial, since scheduling crew already in the tactical phase is preferred over calling upon additional crew last-minute in the operational phase. Even with three scenarios, however, a significant number of excess duties is still required. In other words, it is hard to be prepared for demand peaks in the operational phase that never occur in historical data. Moreover, fixing the capacity plan in the tactical phase comes at a sizeable efficiency loss. In the fully flexible case where all duties can be freely scheduled in the operational phase, only 60.44 duties are required on average. The efficiency loss arises from (i) the need for excess duties and (ii) templates not assigned a duty but left empty. The number of templates left empty can be computed based on Table~\ref{table:effect_scenarios} by adding the number of excess duties to and subtracting the number of total duties from the number of templates. This value is quite stable, ranging around 0.5 for all settings.

\begin{table}[htbp!]
  	\centering 
  	\renewcommand{\arraystretch}{1}
  	\caption{Effect of number of scenarios and maximum number of types.}
  	\label{table:effect_scenarios} 
  	\begin{adjustbox}{max width= \textwidth}
  	\begin{tabular}{@{\extracolsep{0pt}}llrrrr} 
		\toprule 
  		Types & Scenarios & Templates & Duties & Excess & Avg. workload (h) \\
  		\midrule
		 \multirow{3}{*}{10} & 1 & 60.45 & 62.32 & 2.24 & 7.71 \\
		& 2 & 61.82 & 62.40 & 1.13 & 7.70 \\
		& 3 & 62.14 & 62.57 & 1.09 & 7.67 \\
		\midrule
		 \multirow{3}{*}{15} & 1 & 60.18 & 62.14 & 2.50 & 7.72 \\
		& 2 & 61.29 & 61.96 & 1.15 & 7.74 \\
		& 3 & 61.69 & 62.06 & 0.90 & 7.73 \\
		\midrule
		 \multirow{3}{*}{20} & 1 & 60.16 & 62.07 & 2.39 & 7.72 \\
		& 2 & 61.41 & 62.08 & 1.17 & 7.73 \\
		& 3 & 61.91 & 62.09 & 0.88 & 7.72 \\
		\bottomrule 
	\end{tabular} 
  	\end{adjustbox}
\end{table}

Table~\ref{table:effect_scenarios} also shows the cost of obtaining parsimonious solutions. Enforcing a maximum of ten unique template types, as compared to 15, appears to be somewhat restrictive, yielding solutions featuring more chosen templates and duties. Allowing for more than 15 template types does not appear to have a significant benefit. Solutions with up to 20 unique templates are similar in terms of duties and average workload, and actually require slightly more capacity. This results from the slightly higher optimality gap on these relatively unrestricted instances (see Section~\ref{sec:computational}). All in all, a setting with three scenarios and up to 15 template types seems to be preferred, validating this as our choice of default setting in the remainder of the experiments.

\begin{figure}[htbp!]
\centering
 \begin{adjustbox}{max width= \textwidth}
\begin{tikzpicture}
\begin{groupplot}[group style={group name=my plots,group size= 1 by 4,horizontal sep =1.5cm, vertical sep =1cm}, width=18cm, height=7cm, groupplot xlabel={Hour}, groupplot ylabel={Templates (Excess duties)}]
       \nextgroupplot[
	ylabel={Asd, Tuesday},
         ymin =-25,     
	 ymax=75,  
	 xmin=4, 
	xmax=28,
         minor x tick num = 0,
         minor tick length=2ex,
         ytick = {-25, 0, 25, 50, 75},
	legend pos = north east,
	legend cell align = left,
        ]
        \addlegendimage{empty legend}
        \addplot +[style={black, mark=none, thick}] 
table [col sep= semicolon, x index = 0, y index=1] {\subfix{templates_Asd_Tue_1.csv}};

 	\addplot +[style={black, dashed, mark=none, thick}]
table [col sep= semicolon, x index = 0, y index=1] {\subfix{templates_Asd_Tue_2.csv}};

	\addplot +[style={black, dotted, mark=none, thick}]
table [col sep= semicolon, x index = 0, y index=1] {\subfix{templates_Asd_Tue_3.csv}};

	\addplot +[style={black, mark=none, thick}]
table [col sep= semicolon, x index = 0, y index=2]{\subfix{templates_Asd_Tue_1_negative.csv}};

	\addplot +[style={black, dashed, mark=none, thick}]
table [col sep= semicolon, x index = 0, y index=2] {\subfix{templates_Asd_Tue_2_negative.csv}};

	\addplot +[style={black, dotted, mark=none, thick}]
table [col sep= semicolon, x index = 0, y index=2] {\subfix{templates_Asd_Tue_3_negative.csv}};

	\addlegendentry{\hspace{-.6cm}\# Scenarios}
   	\addlegendentry{1}
   	\addlegendentry{2}
	\addlegendentry{3}

\nextgroupplot[
         ymin =-30,     
	   ymax=45,  
	   xmin=4,
	   xmax=28,
	ylabel={Asd, Saturday},
        minor x tick num = 0,
        minor tick length=2ex,
         ytick = {-30, -15, 0, 15, 30, 45},
        ]
        \addplot +[style={black, mark=none, thick}] 
table [col sep= semicolon, x index = 0, y index=1] {\subfix{templates_Asd_Sat_1.csv}};

 	\addplot +[style={black, dashed, mark=none, thick}]
table [col sep= semicolon, x index = 0, y index=1] {\subfix{templates_Asd_Sat_2.csv}};

	\addplot +[style={black, dotted, mark=none, thick}]
table [col sep= semicolon, x index = 0, y index=1] {\subfix{templates_Asd_Sat_3.csv}};

	\addplot +[style={black, mark=none, thick}]
table [col sep= semicolon, x index = 0, y index=2] {\subfix{templates_Asd_Sat_1_negative.csv}};

	\addplot +[style={black, dashed, mark=none, thick}]
table [col sep= semicolon, x index = 0, y index=2] {\subfix{templates_Asd_Sat_2_negative.csv}};

	\addplot +[style={black, dotted, mark=none, thick}]
table [col sep= semicolon, x index = 0, y index=2] {\subfix{templates_Asd_Sat_3_negative.csv}};

\nextgroupplot[
         ymin =-15,     
	ymax=45,  
	xmin=4,
	xmax=28,
	ylabel={Nm, Tuesday},
        minor x tick num = 0,
        minor tick length=2ex,
         ytick = {-15, 0, 15, 30, 45},
	legend pos = north west,
	legend cell align = left,
        ]
        \addlegendimage{empty legend}
        \addplot +[style={black, mark=none, thick}] 
table [col sep= semicolon, x index = 0, y index=1] {\subfix{templates_Nm_Tue_1.csv}};

 	\addplot +[style={black, dashed, mark=none, thick}]
table [col sep= semicolon, x index = 0, y index=1] {\subfix{templates_Nm_Tue_2.csv}};

	\addplot +[style={black, dotted, mark=none, thick}]
table [col sep= semicolon, x index = 0, y index=1] {\subfix{templates_Nm_Tue_3.csv}};

	\addplot +[style={black, mark=none, thick}]
table [col sep= semicolon, x index = 0, y index=3] {\subfix{templates_Nm_Tue_1.csv}};

	\addplot +[style={black, dashed, mark=none, thick}]
table [col sep= semicolon, x index = 0, y index=3] {\subfix{templates_Nm_Tue_2.csv}};

	\addplot +[style={black, dotted, mark=none, thick}]
table [col sep= semicolon, x index = 0, y index=3] {\subfix{templates_Nm_Tue_3.csv}};

\nextgroupplot[
         ymin =-10,     
	ymax=30,  
	xmin=4,
	xmax=28,
	ylabel={Nm, Saturday},
        minor x tick num = 0,
        minor tick length=2ex,
         ytick = {-10, 0, 10, 20, 30},
        ]
        \addplot +[style={black, mark=none, thick}] 
table [col sep= semicolon, x index = 0, y index=1] {\subfix{templates_Nm_Sat_1.csv}};

 	\addplot +[style={black, dashed, mark=none, thick}]
table [col sep= semicolon, x index = 0, y index=1] {\subfix{templates_Nm_Sat_2.csv}};

	\addplot +[style={black, dotted, mark=none, thick}]
table [col sep= semicolon, x index = 0, y index=1] {\subfix{templates_Nm_Sat_3.csv}};

	\addplot +[style={black, mark=none, thick}]
table [col sep= semicolon, x index = 0, y index=3] {\subfix{templates_Nm_Sat_1.csv}};

	\addplot +[style={black, dashed, mark=none, thick}]
table [col sep= semicolon, x index = 0, y index=3] {\subfix{templates_Nm_Sat_2.csv}};

	\addplot +[style={black, dotted, mark=none, thick}]
table [col sep= semicolon, x index = 0, y index=3] {\subfix{templates_Nm_Sat_3.csv}};
\end{groupplot}
\end{tikzpicture}
\end{adjustbox}
\caption{Temporal distribution of templates (positive part of y-axis) and excess duties (negative part of y-axis) for crew bases Amsterdam (Asd) and Nijmegen (Nm) on Tuesdays and Saturdays in the last four weeks of the planning period. At most 15 unique templates are allowed.}
\label{fig:templates_excess}
\end{figure}
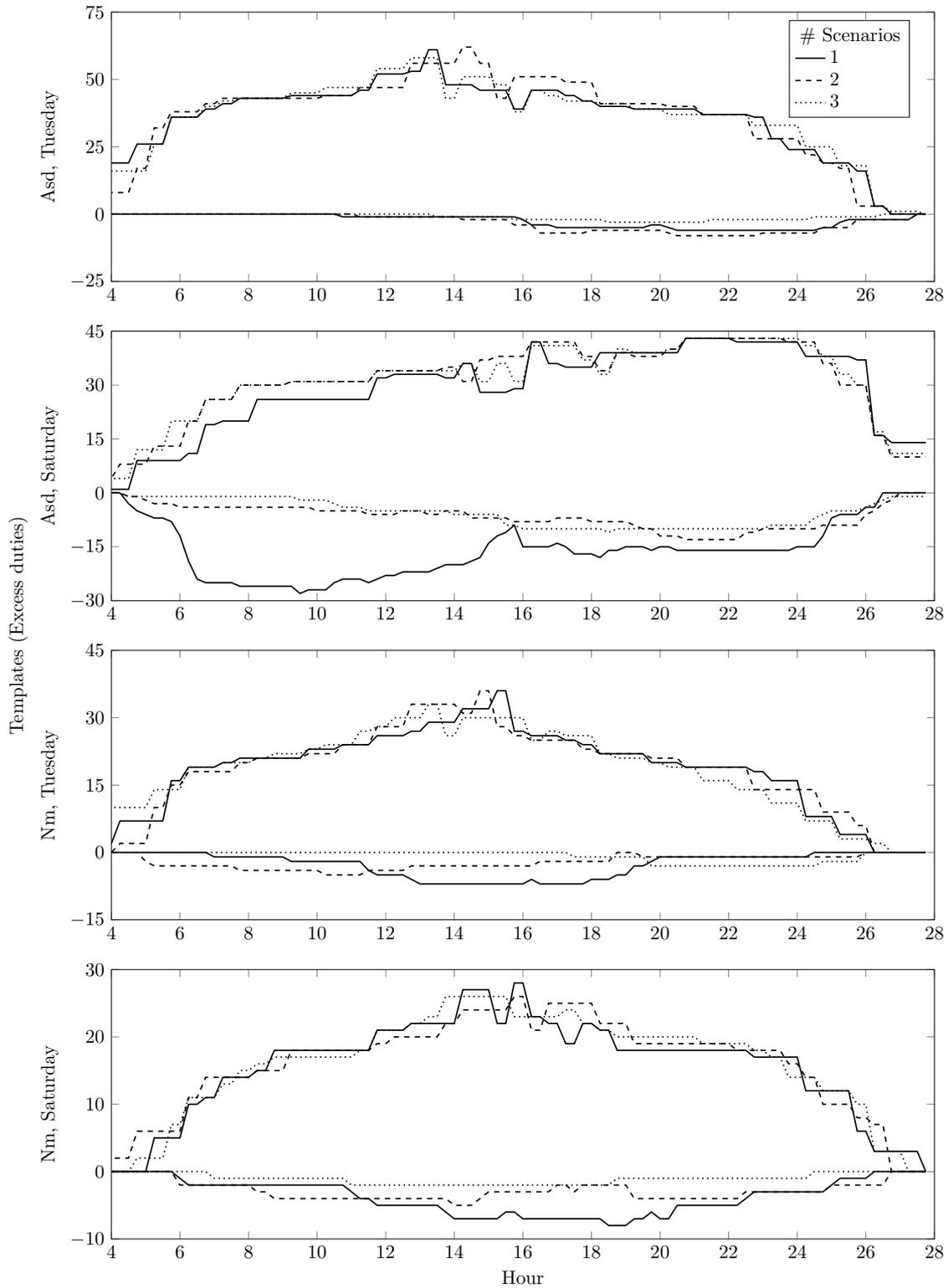

To gain more insight into the effect of the scenario set on solution quality, Figure~\ref{fig:templates_excess} shows the temporal distribution of capacity on the same four instances as Figure~\ref{fig:scenarios} for different numbers of scenarios. For each instance, number of scenarios, and time of the day, it shows the number of active scheduled template (positive part of y-axis) and the number of active excess duties (negative part of y-axis). The latter number is computed by aggregating over all four weeks of the evaluation period. Enlarging the scenario set leads to the scheduling of additional templates at very specific moments of the day, typically early in the morning or late in the evening. It also happens that redundant capacity is shifted away from other times of the day, as illustrated on Tuesday morning in Amsterdam or Tuesday evening in Nijmegen. All in all, adding scenarios leads to a strong reduction in the need for excess duties. Peaks in the number of required excess duties are strongly reduced or even eliminated completely.

\subsection{Reserves}
\label{subsec:practical_reserves}

Table~\ref{table:effect_reserves} presents the effect of varying the percentage of reserve templates. Recall that reserve templates span the whole day, ideally suited for accommodating large demand shocks in the set of tasks to be performed. This is reflected in the results in Table~\ref{table:effect_reserves}, where we find that allowing for some reserves is crucial in attaining efficient schedules. Increasing the percentage of reserves from zero to ten percent leads to a major reduction in the number of required templates, regular duties, and excess duties. The average workload also increases, indicating that reserve templates are assigned efficient duties with a high workload. The marginal gain of adding reserves drops sharply, however, and there appears to be no sizeable benefit in employing more than ten percent reserves. A large fraction of reserve templates is also undesirable from a practical point of view, as it requires crew members to be on stand-by for the whole day. Finally, the results also allow us to estimate the effect of imposing template restrictions. By comparing the first row of Table~\ref{table:effect_reserves} with the last row of Table~\ref{table:benchmark_second}, we find that on average more than one template fewer and fewer excess duties are required when omitting all template restrictions.  

\begin{table}[htbp!]
  	\centering 
  	\renewcommand{\arraystretch}{1}
  	\caption{Effect of reserve templates.}
  	\label{table:effect_reserves} 
  	\begin{adjustbox}{max width= \textwidth}
  	\begin{tabular}{@{\extracolsep{0pt}}lrrrr} 
		\toprule 
  		Reserves (\%) & Templates & Duties & Excess & Avg. workload (h) \\
  		\midrule
		0 & 63.39 & 63.66 & 1.36 & 7.58 \\
		5 & 61.86 & 62.30 & 1.07 & 7.70 \\
		10 & 61.79 & 62.06 & 0.90 & 7.73\\
		15 & 61.79 & 62.15 & 0.89 & 7.71 \\
		20 & 61.71 & 62.05 & 0.86 & 7.72 \\
		\bottomrule 
	\end{tabular} 
  	\end{adjustbox}
\end{table}

\subsection{Template Length}
\label{subsec:practical_length}

Table~\ref{table:effect_length} shows the effect of varying the template length, displaying results for template lengths ranging from eight hours and 45 minutes to ten hours in steps of 15 minutes. When the template length is reduced to nine hours or less, the number of templates and duties grows rapidly. The average workload drops sharply, as the short templates do not allow for efficient, long duties. A reversed, but weaker effect is observed when increasing the template length beyond nine and a half hours. The additional flexibility of longer templates allows for more efficient crew scheduling solutions. Less capacity is required in the tactical phase, and fewer duties are scheduled in the operational phase. Notably, the number of excess duties is not strongly affected by the template length, except in the case of a ten-hour template length. The average workload strongly increases when the template length is increased from nine hours by half an hour or more, indicating that more efficient duties are scheduled and less total workload is required to cover all tasks. From the crew's point of view, however, shorter template lengths are strongly preferred, as they make it easier to schedule their personal lives around the reserved working hours. All in all, a template length of nine and a half hours seems to adequately balance solution quality and crew preferences.

\begin{table}[htbp!]
  	\centering 
  	\renewcommand{\arraystretch}{1}
  	\caption{Effect of template length.}
  	\label{table:effect_length} 
  	\begin{adjustbox}{max width= \textwidth}
  	\begin{tabular}{@{\extracolsep{0pt}}lrrrr} 
		\toprule 
  		Length (h) & Templates & Duties & Excess & Avg. workload (h) \\
  		\midrule
		8:45 & 67.02 & 67.34 & 0.96 & 7.25 \\
		9:00 & 64.11 & 64.49 & 0.89 & 7.49 \\
		9:15 & 62.07 & 62.44 & 0.98 & 7.69 \\
		9:30 & 61.79 & 62.06 & 0.90 & 7.73 \\
		9:45 & 61.75 & 62.03 & 0.85 & 7.72 \\
		10:00 & 61.80 & 61.97 & 0.76 & 7.73 \\
		\bottomrule 
	\end{tabular} 
  	\end{adjustbox}
\end{table}

\subsection{Day of Week}
\label{subsec:practical_day}

Table~\ref{table:effect_day} illustrates the strong heterogeneity in solution quality across days of the week, reflecting the variety observed earlier in Table~\ref{table:instances}. While the lowest number of templates is scheduled in the weekend, the number of excess duties is relatively highest on Saturday and Sunday. This reflects the fact that maintenance works are typically scheduled over the weekend, as observed in the low similarity ratio in Table~\ref{table:instances}. In contrast, few excess duties are required from Tuesday to Thursday, days with a high similarity ratio where demand appears to be more stable and scenarios are more representative of future operations. 

\begin{table}[!htbp]
  	\centering 
  	\renewcommand{\arraystretch}{1}
  	\caption{Effect of day of week.}
  	\label{table:effect_day} 
  	\begin{adjustbox}{max width= \textwidth}
  	\begin{tabular}{@{\extracolsep{0pt}}lrrrr} 
		\toprule 
  		Day & Templates & Duties & Excess & Avg. workload (h) \\
  		\midrule
		Mon & 63.13	& 63.66 & 0.84 & 7.74 \\
		Tue & 63.50	& 63.88 & 0.47 & 7.66  \\ 
		Wed & 62.50	& 62.84 & 0.47 & 7.82  \\
		Thu & 70.13	& 70.34 & 0.38 & 7.73 \\ 
		Fri & 70.13 & 70.69 & 1.28 & 7.68 \\
		Sat & 57.88 & 57.84 & 1.84 & 7.70 \\
		Sun & 45.25	& 45.19 & 1.03 & 7.79 \\
		\bottomrule 
	\end{tabular} 
  	\end{adjustbox}
\end{table}

\FloatBarrier

\section{Conclusion}
\label{sec:conclusion}

Motivated by Netherlands Railways, we consider robust tactical railway crew scheduling for a large passenger railway operator. The goal is to select a cost-efficient set of templates in the tactical planning phase, that is robust with respect to uncertainty about the work to be performed in the operational planning phase. Building on the work of \cite{rahlmann2021robust}, we impose several template rostering constraints, and propose a two-phase accelerated Benders decomposition algorithm that is able to incorporate these restrictions. Experiments on real-life instances of Netherlands Railways show the efficacy of our method and provide practical insight in how to choose robust templates. Historical planning information can be effectively leveraged to obtain robust templates, parsimonious solutions can be obtained at low extra costs, and a high number of excess duties can be avoided by using templates of reasonable length and a minimum number of reserve templates. Moreover, our method compares favourably with a literature benchmark on instances without rostering constraints, solving three times as many large instances to optimality. This comes at the cost of a doubling in average computing time, but computing time is typically not a bottleneck in the tactical planning phase.

Some care must be taken when extrapolating our results to the full network of Netherlands Railways. First, our experiments show that a small number of excess duties is still required on most instances. This necessitates the need for stand-by crew, on top of the crew members already assigned a reserve template. Second, we focus on medium to large crew bases of the Netherlands. Constructing robust rosters at small stations can be more challenging, as minor timetable changes can have a relatively bigger impact. Third, the template rostering constraints do not guarantee that efficient template-based rosters can indeed be constructed using the chosen templates. A simulation of the complete crew planning process is required to determine whether this is the case. 

This work also provides interesting directions for future research. While we have constructed scenarios from historical planning data in a straightforward manner, one can imagine that more sophisticated, data-driven approaches offer more representative scenario sets. Alternatively, one could consider imposing some smoothness conditions on the temporal distribution of template capacity to avoid overfitting on historical data. Moreover, integrating the template selection step with crew rostering, in the line of \cite{borndorfer2017integration}, would be a challenging but promising exercition. Although the method proposed in this work is not yet computationally tractable for the full Dutch railway network, it might lend itself to geographic decomposition strategies or meta-heuristic approaches. Finally, we believe that our method can be applied to other crew planning problems where the construction of a preliminary capacity plan precedes a more detailed work assignment.

\section*{Acknowledgement}

We thank G\'abor Mar\'oti for fruitful discussions of earlier versions of this work.

\bibliographystyle{abbrvnat}
\bibliography{references}

\end{document}